\documentclass[extra]{gji}

\usepackage{timet,color}
\usepackage[urlcolor=blue,citecolor=black,linkcolor=black]{hyperref}
\usepackage{amssymb,amsfonts} % math equations
\usepackage{booktabs} % create professional-looking tables
\usepackage{graphicx} % insert pictures via \includegraphics
\usepackage{algorithm, algorithmicx, algpseudocode} % pseudo codes
\usepackage{layouts}  % Allows displaying layout values
\usepackage[fleqn]{amsmath}
\usepackage{subcaption}

\renewcommand{\vec}[1]{\boldsymbol{#1}}
\newcommand{\src}{s}
\newcommand{\fwd}{\vec{f}}
\newcommand{\data}{\vec{d}}
\newcommand{\postcov}{\Gamma_{\mathrm{post}}}
\newcommand{\prcov}{\Gamma_{\mathrm{pr}}}
\newcommand{\likecov}{\Gamma_{\mathrm{noise}}}
\newcommand{\map}{\vec m_{\mathrm{map}}}
\newcommand{\Hmap}{H_{\mathrm{map}}}

\newcommand{\argmax}{\mathop{\mathrm{arg\,max}}}
\title[Seismic UQ with high-rank Hessian approximations]
  {Accelerating seismic inversion and uncertainty quantification with efficient high-rank Hessian approximations}
\author[Mathew Hu]
  {Mathew Hu$^1$, Nick Alger$^1$, Rami Nammour$^2$, Omar Ghattas$^{1,3}$\\
  $^1$Oden Institute for Computational Engineering and Sciences, The University of Texas at Austin, Austin, TX 78712, USA\\
  $^2$TotalEnergies E\&P USA Inc., 1201 Louisiana Street, Suite 1800, Houston, Texas 77002, USA\\
  $^3$Walker Department of Mechanical Engineering, The University of Texas at Austin, Austin, TX 78712, USA
  }
\date{}
\pagerange{\pageref{firstpage}--\pageref{lastpage}}
\volume{VVV}
\pubyear{2025}

\begin{document}

\label{firstpage}

\maketitle

\begin{summary}
Efficient high-rank approximations of the Hessian (normal operator) can accelerate seismic full-waveform inversion (FWI) and uncertainty quantification (UQ). In FWI, the Hessian encodes the local curvature of the FWI objective function. Approximations of the inverse of the Hessian may be used as a preconditioner for inexact Newton-Krylov algorithms or quasi-Newton algorithms, significantly reducing computational costs and improving recovery in deeper subsurface regions. In Bayesian UQ, the inverse of the Hessian is the covariance of a local Gaussian approximation to the posterior. Hessian approximations enable the construction of Markov chain Monte Carlo (MCMC) proposal distributions that faithfully capture the directional scalings of the posterior, thereby enhancing the efficiency of MCMC algorithms. Computing the exact Hessian is intractable for large-scale problems because the Hessian is accessible only matrix-free, i.e. through matrix-vector products, and performing each matrix-vector product requires costly solution of wave equations. Moreover, the Hessian is high-rank, which means that low-rank Hessian approximation methods, often employed to make large-scale inverse problems tractable, are inefficient. We adapt two existing high-rank Hessian approximations---the point spread function (PSF) method and the pseudo-differential operator probing (PDO) method. Building on an observed duality between these two approaches, we further develop a novel approximation method that unifies their complementary strengths. We validate these methods on a synthetic quadratic model and on the Marmousi model. Numerical experiments demonstrate that these high-rank Hessian approximations substantially reduce the computational costs in FWI. In UQ, we observe that high-rank Hessian approximations are essential. MCMC samples computed using no Hessian approximation or a low-rank approximation explore the posterior extremely slowly, providing little meaningful statistical information after tens of thousands of iterations and causing the variance to be underestimated. At the same time, the effective sample size (ESS) is overestimated, providing false confidence. In contrast, MCMC samples generated using the high-rank Hessian approximations provide meaningful statistical information about the posterior and more accurately assess the posterior variance.
\end{summary}


\section{Introduction}

Seismic inversion plays a pivotal role in numerous applications, including subsurface exploration \citep{sheriff1995exploration,barclay2008seismic}, reservoir management \citep{lumley2001time}, carbon sequestration \citep{um2023real}, and scientific investigation of the Earth’s structure.
Traditional deterministic seismic inversion methods \citep{lailly1983seismic,tarantola1984inversion,virieux2009overview}  yield a single ``best fit'' model of the subsurface that minimizes the misfit between observed and predicted seismic data.
However, seismic inversion is inherently ill-posed in the sense that many distinct models can generate predicted observations that fit the observed data to within the noise level. Consequently, it is desirable to not only compute an estimate of the unknown model parameters, but also to quantify the uncertainty in the estimate in order to assess the reliability and variability of inversion results. This information is critical for evaluating risk and informing decision-making.

Bayesian inference offers a promising framework for uncertainty quantification (UQ) in seismic inversion \citep{martin2012stochastic,stuart2010inverse,tarantola2005inverse}. In this approach, the posterior probability distribution over all possible subsurface models is given by Bayes' rule, and statistical information is extracted from the posterior using methods such as Markov Chain Monte Carlo (MCMC) sampling. 
However, existing MCMC sampling methods often suffer from inefficiency due to the high-dimensional and ill-conditioned nature of the posterior distribution. In such cases, small sampling steps must be taken to maintain reasonable acceptance rates, resulting in poor exploration and requiring a large number of samples to obtain accurate statistical estimates. Since each sample involves evaluating the posterior, MCMC becomes computationally expensive.

Fortunately, the number of samples required for effective posterior exploration can be significantly reduced by incorporating Hessian information into the proposal distribution. 
Hessian-informed proposals, such as those based on the Laplace approximation, align proposal distributions with the local geometry of the posterior, enabling larger and more targeted steps. A variety of algorithms have been developed to exploit this idea. For example, the Hessian-preconditioned MALA (H-MALA) utilizes the curvature information to capture the posterior correlation \citep{martin2012stochastic}. The dimension-independent likelihood-informed (DILI) method combines likelihood-informed subspace techniques with MCMC to construct efficient proposals informed by the Hessian \citep{cui2016dimension}. A recent comprehensive assessment of Hessian-informed MCMC methods is given in \citet{kim2023hippylib}, which evaluates a broad spectrum of proposals using low-rank Hessian approximations and shows substantial speedups across PDE-constrained inverse problems. To ensure performance that is independent of discretization resolution, it is important to combine Hessian-informed proposals with dimension-invariant sampling strategies. For this, we adopt the generalized preconditioned Crank-Nicolson (gpCN) algorithm \citep{pinski2015algorithms,rudolf2018generalization,cao2024efficient}, which generalizes the classical pCN method \citep{cotter2013mcmc} by allowing proposals from an arbitrary Gaussian distribution. Other dimension-invariant Hessian-informed samplers include the DILI \citep{cui2016dimension} and the Hessian-informed infinite-dimensional MALA (H-$\infty$-MALA) \citep{beskos2017geometric,kim2023hippylib}.

However, computing the full Hessian matrix is intractable for large-scale problems due to both computational and memory limitations. Therefore, efficient Hessian approximations are essential.
Low-rank approximation methods have demonstrated success in various large-scale inverse problems and UQ tasks \citep{isaac2015scalable,petra2014computational,bui2013computational,martin2012stochastic}. Yet, in seismic settings, data becomes highly informative when a large amount of sources and receivers are deployed. Highly informative data leads to a high-rank Hessian, which renders low-rank approaches inefficient \citep{alger2017data,alger2019data}. Alternative methods that approximate the diagonal of the Hessian or its inverse \citep{beylkin1985imaging,beydoun1989elastic,yang2018time} can yield high-rank approximations, but often neglect the spatial spreading effects of the operator. Other high-rank approaches have exploited the sparsity of the Hessian in curvelet frames \citep{herrmann2008sparsity,herrmann2009curvelet}, the locality and patterns of the point spread function \citep{tang2009target,zhu2016bayesian,alger2019scalable,yang2021approximating,yang2022efficient,alger2024point}, the pseudo-differential properties of the Hessian's symbol \citep{bao1996computation,nammour2009approximate,nammour2011approximate,demanet2012matrix}, or the Hessian's hierarchical structure \citep{ambartsumyan2020hierarchical,hartland2023hierarchical}. More recently, machine learning techniques have shown promise in approximating the seismic Hessian \citep{alfarhan2024robust,kaur2020improving}. 

Recent advances in Hessian approximation have primarily targeted the acceleration of deterministic seismic inversion. These techniques aim to address the common challenge that deeper regions of the subsurface model are less informed by the data and therefore harder to recover.
In contrast, the use of Hessian approximations for UQ remains relatively unexplored. Although both deterministic inversion and UQ rely on Hessian information, the requirements for UQ are substantially more stringent.
First, deterministic inversion only requires a preconditioner, typically an approximation of the inverse Hessian, which does not need to be symmetric positive definite (SPD). For UQ, however, the Hessian approximation itself must be SPD, and its inverse square root must be accurately approximated to enable Gaussian sampling from the Laplace approximation of the posterior.
Second, in deterministic inversion, the preconditioner is typically applied only to the gradient of the objective function, which lies within the range of the migration operator $F^*$, as shown in Equation~(\ref{equ:grad}). Thus, the preconditioner only needs to approximate the curvature along a restricted set of directions. In contrast, UQ involves applying the inverse square root of the Hessian approximation to random vectors, requiring accurate curvature information across all directions.
Consequently, Hessian approximations that are effective for deterministic inversion may not be suitable for uncertainty quantification.

In this work, we accelerate seismic uncertainty quantification by adapting existing high-rank Hessian approximation techniques and introducing a novel algorithm. The proposed algorithm is based on a synergistic combination of the point spread function method and the pseudo-differential operator probing method, leveraging a hidden duality between these seemingly distinct approaches. To validate the effectiveness of the accelerated UQ framework, we first apply it to an idealized quadratic model, where the true uncertainty is analytically known. This setting allows us to clearly demonstrate the benefits of high-rank Hessian approximations. We then extend the approach to the more realistic synthetic Marmousi model to assess its performance in complex scenarios. Our results show that high-rank Hessian approximation is essential in MCMC-based UQ---not only for improving convergence rates but also for ensuring the reliability of uncertainty estimates. Without such high-rank Hessian approximation, MCMC samples may appear to converge, yet still produce significantly underestimated variances. This misrepresentation of uncertainty can lead to overconfident and potentially misleading conclusions, with serious implications for risk assessment and decision-making in practice.

\section{Seismic Inversion Bayesian Framework}

Assuming constant density and isotropic bulk modulus, the governing acoustic wave equation for the seismic inversion is given by
\begin{equation} \label{equ:pde}
    \frac{1}{m^2(\vec x)} \frac{\partial^2 u(\vec x, t)}{\partial^2 t} - \nabla^2 u(\vec x, t) = \src(t) \delta(\vec x - \vec x_s),
\end{equation}
where $\vec{x}$ and $t$ denote spatial and temporal variables, respectively. Here, $m(\vec{x})$ represents the wave speed (the unknown model parameter) and $u(\vec{x}, t)$ is the resulting pressure field. The right-hand side $s(t)$ denotes the source function, and $\vec x_s$ is the source location. 

Given a source location $\vec x_s$ and source function $s(t)$, and assuming a model $m(\vec x)$, one can numerically solve the wave equation, for instance using finite difference methods, to compute the pressure field $u(\vec{x}, t)$. The observed data are recorded as values of $u(\vec{x}, t)$ at receiver locations. This process defines a forward map $\vec f$ that maps each discretized model $\vec m$ to observed data $\data$. Assuming additive Gaussian noise, the noisy data is modeled as
\begin{equation}\label{equ:data}
    \data = \vec f(\vec m) + \vec \eta, \quad \vec \eta \sim \mathcal{N}(0, \Gamma_{\mathrm{noise}}),
\end{equation}
where $\vec \eta$ is a Gaussian random variable with zero mean and covariance $\Gamma_{\mathrm{noise}}$. In seismic inversion, the data $\data$ is given, and the goal is to recover the model $\vec m$ and quantify its uncertainty.

Bayes' rule combines a prior distribution $\pi_{\mathrm{pr}}(\vec m)$ with a likelihood function $\pi_{\mathrm{like}}(\data|\vec m)$ to form the posterior distribution over all possible models:
\begin{equation}
    \pi_{\mathrm{post}}(\vec m|\data) = \frac{1}{Z} \pi_{\mathrm{like}}(\data | \vec m) \pi_{\mathrm{pr}}(\vec m),
\end{equation}
where the normalization constant (also called the evidence) $Z$ is defined as
\begin{equation}
    Z = \int \pi_{\mathrm{like}}(\data|\vec m) \pi_{\mathrm{pr}}(\vec m) d \vec m.
\end{equation}

The likelihood $\pi_{\mathrm{like}}(\data|\vec m)$ quantifies the probability of observing the data $\data$ given the model parameters $\vec m$. Under the additive Gaussian noise assumption~(\ref{equ:data}), the likelihood yields
\begin{equation}
    \pi_{\mathrm{like}}(\data|\vec m) \propto 
    \exp\left(
    -\frac12 \| \vec f (\vec m)-\data\|^2_{\likecov^{-1}}
    \right),
\end{equation}
where the weighted norm $\|\cdot\|_A$ is defined as $\|\vec v\|_A = \sqrt{\vec v^T A \vec v}$.

The prior encodes knowledge about the model before observing any data. We adopt a Mat\'{e}rn class prior, which is widely used in physical modeling for controlling both marginal variance and smoothness. This prior is defined as a Gaussian distribution:
\begin{equation}
    \pi_{\mathrm{pr}}(\vec m) \propto
    \exp\left(
    -\frac12\|\vec m-\vec m_{\mathrm{pr}}\|^2_{\prcov^{-1}}
    \right),
\end{equation}
where $\vec m_{\mathrm{pr}}$ is the prior mean, and $\prcov$ is a Mat\'{e}rn covariance operator \citep{lindgren2011explicit}. This covariance is equivalent (up to the boundary conditions) to the discretization of the inverse of the $\alpha$-th power of an elliptic differential operator:
\begin{equation}
    \Gamma_{\mathrm{pr}} = (\delta I_h - \gamma \Delta_h)^{-\alpha},
\end{equation}
where $\alpha=2$ is chosen to ensure that the prior covariance is a trace-class operator in two or three dimensions \citep{stuart2010inverse}, resulting in the so-called biharmonic prior. The positive constants $\delta$ and $\gamma$ control the correlation length $\rho \propto \sqrt{\gamma/\delta}$ 
 (the distance for which the two-point correlation is 0.1) and the prior’s pointwise variance $\sigma \propto 1 / \sqrt{\delta \gamma}$ \citep{daon2018mitigating}.

Combining the likelihood and prior through Bayes' rule yields the posterior distribution:
\begin{align}
    &\pi_{\mathrm{post}}(\vec m|\data) \propto \nonumber\\
    &\exp\left(
    -\frac12 \| \vec f (\vec m)-\data\|^2_{\likecov^{-1}}
    -\frac12\|\vec m - \vec m_{\mathrm{pr}}\|^2_{\prcov^{-1}}
    \right).
\end{align}

\subsection{Maximum a posteriori}

Given the posterior distribution, a key quantity of interest is the maximum a posteriori (MAP) point, denoted by $\vec m_{\mathrm{map}}$, which maximizes the posterior:
\begin{equation}
    \map:= \argmax_{\vec m} ~\pi_{\mathrm{post}}(\vec m|\data).
\end{equation}
Equivalently, $\vec m_{\mathrm{map}}$ minimizes the negative log-posterior. Thus, it can be obtained by solving the following optimization problem:
\begin{equation}
    \min_{\vec{m}} \Phi(\vec m) := \frac12 \| \fwd (\vec m)-\data\|^2_{\likecov^{-1}}+\frac12\|\vec m-\vec m_{\mathrm{pr}}\|^2_{\prcov^{-1}}.
\end{equation}
This MAP point $\map$ is equivalent to the ``best-fit'' model in the context of full-waveform inversion (FWI). Specifically, the inverse prior covariance acts as a regularization operator, while the inverse noise covariance serves as a weighting for the $L^2$ data misfit term when $\likecov = \sigma^2 I$.

This optimization problem can be solved using iterative algorithms, provided that the gradient is available. The gradient of the objective function is given by:
\begin{equation}\label{equ:grad}
    \nabla \Phi(\vec m) = F^* \likecov^{-1} (\fwd(\vec m)-\data ) + \prcov^{-1}(\vec m-\vec m_{\mathrm{pr}})
\end{equation}
where $F=\nabla \fwd$ denotes the Born modeling operator. The adjoint operator $F^*$, commonly referred to as the migration operator, maps data residuals back into the model space.

We employ the L-BFGS algorithm \cite[Section 7.2]{nocedal1999numerical} to solve this optimization problem and further enhance its performance by initializing the internal inverse Hessian approximation using our proposed Hessian approximations. This preconditioning strategy improves the inversion; for example, it enables more balanced recovery of both shallow and deep regions of the model, effectively addressing the challenge of reconstructing the full subsurface structure.

\subsection{Laplace approximation}

The Laplace approximation to the posterior distribution is the Gaussian distribution $\mathcal{N}(\map, \postcov)$ centered at the MAP point $\map$, with covariance $\postcov = (\nabla^2\Phi(\map))^{-1}$:
\begin{equation}
    \pi_{\mathrm{post}}^{\mathcal{L}}(\vec m|\data) \propto
    \exp\left( -\frac{1}{2}\|\vec m-\vec m_{\mathrm{map}}\|^2_{\postcov^{-1}} \right).
\end{equation}
This approximation arises from a second-order Taylor expansion of the negative log-posterior about the MAP point, as discussed in Ghattas \& Willcox \shortcite[Section 4.3]{ghattas2021learning}. Since the exact Hessian $\nabla^2\Phi(\map)$ is typically unavailable or computationally expensive to evaluate, we instead use an approximation $\Hmap \approx \nabla^2\Phi(\map)$, yielding an approximate posterior covariance $\widetilde{\Gamma}_{\mathrm{post}} = \Hmap^{-1}$. The resulting approximate posterior distribution is then given by
\begin{equation}\label{equ:la-post}
    \widetilde{\pi}_{\mathrm{post}}^{\mathcal{L}}(\vec m|\data) \propto
    \exp\left( -\frac{1}{2}\|\vec m- \map \|^2_{\widetilde{\Gamma}_{\mathrm{post}}^{-1}} \right).
\end{equation}

\subsection{MCMC-gpCN}

The generalized preconditioned Crank-Nicolson (gpCN) algorithm leverages a Gaussian approximation of the posterior to generate efficient proposals for MCMC sampling. At each iteration $k$, a new proposal $\vec m'$ is generated by perturbing the current state $\vec m_k$ using a sample from the Laplace-approximated posterior: 
\begin{equation}\label{equ:gpcn}
    \vec m' = \map + \sqrt{1-\beta^2}(\vec m_k - \map) + \beta \vec \zeta,
\end{equation}
where $\vec \zeta \sim \mathcal{N}(0, \widetilde{\Gamma}_{\mathrm{post}})$ and $\beta\in(0,1)$ is a tunable step-size.
The proposal is accepted with probability
\begin{equation}
    a(\vec m_k, \vec m') = \min \{1, \exp(\Delta(\vec m_k) - \Delta(\vec m'))\},
\end{equation}
where
\begin{equation}
    \Delta(\vec m) = \Phi(\vec m) - \frac12 \|\vec m - \map\|^2_{\widetilde{\Gamma}_{\mathrm{post}}^{-1}}
\end{equation}
denotes the difference between the negative log posterior and the approximate negative log posterior based on the Laplace approximation (up to a constant shift).
If the Laplace approximation $\widetilde{\pi}_{\mathrm{post}}^{\mathcal{L}}$ closely matches the true posterior, $\Delta(\vec m)$ remains approximately constant, leading to high acceptance rates.
In summary, the gpCN algorithm proceeds as follows:
\begin{enumerate}
    \item Initialize: Set $k=0$ and pick initial sample $\vec m_0$.
    \item Proposal: Generate $\vec m'$ using Equation (\ref{equ:gpcn}).
    \item Accept/Reject: Set $\vec m_{k+1}=\vec m'$ with probability $a(\vec m_k,\vec m')$, or set $\vec m_{k+1}=\vec m_k$ otherwise.
    \item Iterate: $k \leftarrow k+1$ and return to step (ii).
\end{enumerate}

\section{Hessian approximation methods}

The Hessian operator of the objective function $\Phi(\vec m)$ consists of the misfit Hessian $H_d(\vec m)$ and the regularization Hessian $R=\prcov^{-1}$:
\begin{equation}
    \nabla^2\Phi(\vec m) = H(\vec m) = H_d(\vec m) + R
\end{equation}
The regularization Hessian, $\prcov^{-1}$, is either a scaled mass matrix in the case where $L^2$ regularization is used, or a discretization of a differential operator in the case of a Matern prior. In both cases, $\prcov^{-1}$ is a sparse matrix that is known and easy to deal with using conventional approaches. We therefore focus on approximating the misfit Hessian $H_d(\vec m)$, with the understanding that the approximation of $H_d$ is later combined with $R$ to form an approximation of $H$.

Although the full matrix $H_d$ cannot be computed explicitly, it can be accessed via matrix-vector products using adjoint methods \cite[Section 2.6.4]{epanomeritakis2008newton,martin2015computational}.
Each matrix-vector product requires solving two wave equations per source. Thus, computing the full matrix $H_d$ would require $2NN_s$ wave equation solves, where $N$ is the discretization dimension and $N_s$ is the number of sources, resulting in a significant computational burden. Moreover, storing the full matrix $H_d$, which has size $N \times N$, demands substantial memory resources. These challenges highlight the difficulties we face when building Hessian approximations. For Hessian approximation methods to be computationally tractable, they must (1) construct the approximation using information from a small number of matrix-vector products of $H_d$, and (2) yield a Hessian approximation that is represented in a compressed format.

In this section, we introduce low-rank approximation as both a baseline and an auxiliary technique to enhance UQ performance. We then adapt the point spread function method and the pseudo-differential operator probing method, two established Hessian approximation methods that are are widely used in the broader literature. We identify a duality between these methods: the point spread function method can be interpreted as constructing a pseudo-differential operator in a transposed form. Building on this insight, we propose a novel low-rank symbol algorithm that unifies these approaches.

\subsection{Low-rank approximation}

The low-rank approximation method has proven effective in many applications. This method approximates $H$ by forming a truncated eigenvalue decomposition of $H_d$ preconditioned by $R$. Specifically, it approximates the Hessian by solving the generalized eigenvalue problem
\begin{equation}\label{equ:eig}
    H_d v_i= \lambda_i R v_i, \quad \lambda_1 \ge \lambda_2 \ge \cdots \ge \lambda_r,
\end{equation}
for the first $r$ eigenvalues and eigenvectors. Define
\begin{align}
    &V_r = [v_1, v_2, \cdots, v_r],\\
    &\Lambda_r = \mathrm{diag}([\lambda_1, \lambda_2, \cdots, \lambda_r]),
\end{align}
where $V_r^T R V_r = I$. The resulting Hessian approximation and its inverse are given by
\begin{equation}
\widetilde H_d+R = R + R V_r \Lambda_r V_r^T R + \mathcal{O}\left(\sum_{i>r} \lambda_i\right), 
\end{equation}
\begin{equation}
(\widetilde H_d+R)^{-1} = R^{-1} - V_r D_r V_r^T + \mathcal{O}\left(\sum_{i>r} \frac{\lambda_i}{\lambda_i+1}\right),
\end{equation}
where $D_r = \mathrm{diag}([\frac{\lambda_1}{\lambda_1+1}, \frac{\lambda_2}{\lambda_2+1}, \cdots, \frac{\lambda_r}{\lambda_r+1}])$.
Furthermore, one can draw proposal samples $\vec \zeta$ from the Laplace approximation $\mathcal{N}(0, (\widetilde H_d + R)^{-1})$ via
\begin{equation}
    \vec \zeta = (I - V_r S_r V_r^T R) \vec z,
\end{equation}
where $\vec z \sim \mathcal{N}(0, R^{-1})$ and $S_r = I - (\Lambda_r+I)^{-\frac12}$.
The required eigenvalues and eigenvectors are computed via matrix-free methods such as Lanczos \cite[Section 6.3]{saad2011numerical} or randomized SVD \citep{saibaba2016randomized,halko2011finding}. 
This approach is further detailed in the hIPPYlib software framework \citep{villa2021hippylib}, where a randomized SVD algorithm is implemented to solve the generalized eigenvalue problem~(\ref{equ:eig}).

Using these matrix-free methods, the required number of matrix-vector products with $H_d$ scales with the rank of the low-rank approximation. Hence, the low-rank approximation method performs well when the eigenvalues decay quickly. However, in seismic inversion, the eigenvalues decay more slowly with increasing numbers of sources and receivers, and with higher frequencies, making low-rank approximations inefficient. In such cases, a high-rank Hessian approximation method is required.

\subsection{Point spread function method}

The point spread function (PSF) is defined as application of $H_d$ to a delta function.
PSFs are generally localized, enabling the computation of multiple PSFs by applying $H_d$ to a single vector, which consists of a sum of well-separated delta functions. In seismic problems, reflection effects result in local PSFs, whereas transmission effects introduce non-local PSFs. Fortunately, transmission effects are primarily shallow and associated with low frequencies. To mitigate non-locality, we apply a high-pass filter to remove low-frequency components before the approximation, as further discussed in Section 3.7. The locality of the computed PSFs is evident in the implementations shown in Figs \ref{fig:toy-app-psf-vec} and \ref{fig:mar-app-psf-vec}, where PSFs are computed at regularly aligned points, a natural choice in seismic applications \citep{yang2021approximating,yang2022efficient}. 
Let $\mathcal{H}_d: L^2(\mathbb{R}^2) \rightarrow L^2(\mathbb{R}^2)$ denote the undiscretized misfit Hessian operator. The PSF centered at sample point $\vec x_k$ is defined as
\begin{equation}
    p_k(\vec x) = \mathcal{H}_d \delta_{\vec x_k}(\vec x_k + \vec x), \quad k = 1,2,\cdots,r.
\end{equation}
These PSFs are conventionally used to construct a product-convolution approximation of $\mathcal{H}_d$, where a set of interpolation weights, such as radial basis functions, $w_k(\vec x)$ are employed to blend contributions from different convolutions:
\begin{equation}\label{equ:psf-pc}
    \mathcal{H}_d v \approx \sum_{k=1}^r p_k * (w_k \cdot v),
\end{equation}
where ``$*$'' denotes convolution operation and ``$\cdot$'' indicates pointwise multiplication. For a more detailed general discussion of point spread functions and product convolution approximations, we recommend the following papers \citep{escande2017approximation,denis2015fast,gentile2013interpolating}.

\subsection{Pseudo-differential operators}

It has been shown that, under certain conditions---including a smooth background, the exclusion of diving waves, and the dominance of propagating phases---$\mathcal{H}_d$ can be well approximated by a pseudo-differential operator \citep{beylkin1985imaging,stolk2000microlocal,nammour2011multiparameter}.

Pseudo-differential operators ($\Psi$DOs) \citep{hormander2007analysis,demanet2011discrete} are a generalization of differential operators. The operator $\mathcal{H}_d: L^2(\mathbb{R}^2) \rightarrow L^2(\mathbb{R}^2)$ is a $\Psi$DO if it may be written in the form
\begin{equation} \label{equ:psido}
    \mathcal{H}_d v(\vec x) = \int e^{i\vec x \cdot \vec \xi} s(\vec x, \vec \xi) \widehat{v}(\vec \xi) d \vec \xi,
\end{equation}
where $\vec x$ represents spatial coordinates, $\vec\xi$ denotes frequency, $\widehat{v}(\vec \xi)$ is the Fourier transform\footnote{
    Our conventions: 
    $\widehat{v} = \frac{1}{(2\pi)^2} \int e^{-i \vec x\cdot\vec \xi} v(\vec x) d \vec x$, 
    $v = \int e^{i \vec x\cdot\vec \xi} \widehat{v}(\vec \xi) d\vec \xi$.
} of $v$, and the \textit{symbol} $s(\vec x,\vec \xi): \mathbb{R}^2\times\mathbb{R}^2 \to \mathbb{C}$, which defines the scaling in phase space (spatial $\times$ frequency), satisfies certain constraints on its asymptotic behavior in $|\vec \xi|$. Conceptually, a $\Psi$DO operates in three steps: (1) it transforms spatial space to phase space, (2) it applies scaling in phase space via the symbol, and (3) it maps the result back to spatial space. 

The symbol $s(\vec x, \vec \xi)$ is fundamental for the scaling and must belong to a specific symbol class for the $\Psi$DO to be well-defined. For instance, a symbol of order 1, denoted $s \in S^1$, is infinitely differentiable, and all its derivatives are bounded by polynomials of $\vec \xi$ of corresponding orders, i.e.,
\begin{equation}\label{equ:symbol-smooth}
    \left|\partial_{\vec \xi}^\alpha \partial_{\vec x}^\beta s(\vec x, \vec \xi)\right| 
    \leq C_{\alpha\beta} (1+|\vec \xi|)^{1-|\alpha|},
\end{equation}
for all $\vec x, \vec \xi \in \mathbb{R}^2$, multi-indices $\alpha, \beta$, and some constants $C_{\alpha\beta}$. 

Computing all entries of the symbol is computationally prohibitive, and even with a fully known symbol, directly applying the $\Psi$DO to a vector incurs a computational cost of $\mathcal{O}(N^2)$, which exceeds the cost of the wave simulation itself. However, the symbol can often be approximated as low-rank even if the operator corresponding to the symbol (here, $\mathcal{H}_d$) is high-rank. An easy example is the identity operator: it is full rank, but its symbol is rank-1 with all ones. 

The smoothness condition (\ref{equ:symbol-smooth}) enables a separation of the symbol function, which is equivalent to a low-rank approximation of the symbol in terms of spatial vectors and frequency dual vectors:
\begin{equation} \label{equ:lr_symbol}
    s(\vec x,\vec \xi) \approx \sum_{k=1}^r a_k(\vec x) b_k(\vec \xi),
\end{equation}
or equivalently, its discretized matrix form,
\begin{equation} \label{equ:lr_symbol_matrix}
    S \approx AB, \quad A \in \mathbb{R}^{N\times r}, B \in \mathbb{R}^{r\times N}.
\end{equation}
The $\Psi$DO formulation (\ref{equ:psido}) simplifies to
\begin{equation}\label{equ:lr-psido}
    \mathcal{H}_d \approx \sum_{k=1}^r \mathcal{D}_{a_k} \mathcal{F}^{-1} \mathcal{D}_{b_k} \mathcal{F},
\end{equation}
where we use $\mathcal{D}_f$ to denote a diagonal operator that multiplies by the function $f$. From this formula, we see that applying a $\Psi$DO to a vector requires one Fourier transform, $r$ inverse Fourier transforms, and $2r$ pointwise multiplications. The total computational cost of these operations is $\mathcal{O}(rN\log N)$, which makes the pseudo-differential operator method suitable for large-scale applications. We note that \citet{bao1996computation} and \citet{demanet2011discrete} choose different basis functions $a_k(x)$ and $b_k(x)$ but essentially employ the same separation step to derive an efficient algorithm.

In following discussions, we define a \textit{symbol row} as a slice $s(\vec x, \cdot)$, where the spatial variable $\vec x$ is fixed, and a \textit{symbol column} as $s(\cdot, \vec\xi)$, where the frequency variable $\vec \xi$ is fixed. This terminology aligns with the standard definitions of matrix rows and columns. Since both the spatial and frequency space are 2D, the symbol $s(\vec x,\vec \xi)$ is a 4D function. Each symbol row or column is thus a 2D slice of this 4D function, as illustrated in Fig.~\ref{fig:symbol}.

\subsection{Pseudo-differential operator probing method}\label{sec:pdo}

The pseudo-differential operator probing method (PDO) computes the symbol columns $s(\cdot, \vec \xi_k)$ at selected frequency points $\vec \xi_k$ by applying $\mathcal{H}_d$ to sinusoids. Specifically, consider the function $v$ to be the Fourier basis function $\phi_{\vec \xi}(\vec x)=e^{i\vec x \cdot \vec \xi}$ in the definition~(\ref{equ:psido}). Noting that $\widehat{\phi}_{\vec\xi} = \delta_{\vec\xi}$, we obtain
\begin{align}
    &\mathcal{H}_d\phi_{\vec\xi}(\vec x)
    = \int e^{i\vec x \cdot \vec \xi'} s(\vec x, \vec \xi') \delta_{\vec\xi}(\vec \xi') d \vec \xi' \nonumber\\
    &= e^{i\vec x \cdot \vec \xi} s(\vec x, \vec \xi) 
    = \phi_{\vec\xi}(\vec x) s(\vec x, \vec \xi).
\end{align}
Therefore, the symbol columns are expressed as
\begin{equation} \label{equ:symbol}
    s(\cdot,\vec \xi) = \mathcal{H}_d \phi_{\vec \xi} / \phi_{\vec \xi},
\end{equation}
where ``$/$'' denotes pointwise division.
We propose a probing method that computes multiple symbol columns using only one operator application of $\mathcal{H}_d$ by exploiting the locality of $\mathcal{H}_d$ in frequency space as pseudo-differential operators preserve singularities in phase space. %[RN note: this work was never published so this cannot reference my previous work]
Consider a probing vector composed of a sum of sinusoids: 
\begin{equation}\label{equ:pdo-vp}
    v_p(\vec x) = \sum_{k=1}^r \phi_{\vec \xi_k} (\vec x) = \sum_{k=1}^r e^{i\vec x \cdot \vec \xi_k}.
\end{equation}
The application of $\mathcal{H}_d$ to this probing vector turns out to be again a combination of these sinusoids, but each of them is multiplied by a symbol column:
\begin{equation}
    \mathcal{H}_d v_p
    = \sum_{k=1}^r \mathcal{H}_d \phi_{\vec \xi_k}
    = \sum_{k=1}^r \phi_{\vec \xi_k} s(\cdot,\vec \xi_k).
\end{equation}
By applying the Fourier transform to this result $\mathcal{H}_d v_p$, the spectrum of the symbol columns can be found clustered around each frequency sample point $\vec \xi_k$, and hence separated from each other. One can thus extract many symbol columns $s(\cdot,\vec \xi_k), k=1,\cdots, r$, from $\mathcal{H}_d v_p$. 
The locality property can be observed from our implementation in Figs \ref{fig:toy-app-pdo-vec} and \ref{fig:mar-app-pdo-vec}.

Given the computed symbol columns, the PDO method approximates the symbol with a special interpolation strategy that leverages existing knowledge about $\mathcal{H}_d$ as a $\Psi$DO. It has been shown that $\mathcal{H}_d$ is a $\Psi$DO of order 1 in 2D, i.e., $\mathcal{H}_d \in S^1$, under certain conditions \citep{bao1996computation}. 
Based on this asymptotic behavior, the symbol has a polynomial expansion on $\vec \xi$'s magnitude when we consider it in a polar coordinate $\vec \xi = (\rho,\theta)$:
\begin{equation}
    s(\vec x,\vec \xi) = \rho s_1(\vec x,\theta) + s_0(\vec x,\theta) + \rho^{-1} s_{-1}(\vec x,\theta) + \cdots
\end{equation}
Focusing on the principal term leads to an interpolation strategy that extends the symbol linearly in magnitude and uses Fourier interpolation over the angles. That is, for sample frequencies $\vec \xi_k = (\rho_0, \theta_k)$ on a circle of radius $\rho_0$, the interpolation weighting functions are defined as follows:
\begin{equation} \label{equ:pdo-w}
    w_k(\vec \xi) = w_k(\rho, \theta) = \frac{\rho}{\rho_0} w^{\mathrm{circle}}_k(\theta),
\end{equation}
where $w^{\mathrm{circle}}_k(\theta)$ are the angular FFT interpolation weighting functions.

Given symbol columns $s(\vec x,\vec \xi_k)$ and interpolation weighting functions $w_k(\vec \xi)$ for $k=1,2,\cdots,r$, the symbol can be approximated as follows:
\begin{equation}\label{equ:pdo-symbol}
    s(\vec x, \vec \xi) \approx \sum_{k=1}^r s(\vec x,\vec \xi_k) w_k(\vec \xi),
\end{equation}
which aligns with the low-rank symbol in expression~(\ref{equ:lr_symbol}), where the symbol columns are $a_k$ and the weighting functions are $b_k$. Therefore, the resulting approximated misfit Hessian can be applied efficiently through chains of operators illustrated in expression~(\ref{equ:lr-psido}). Additional details of the PDO method are provided in the Appendix.

\subsection{PSF method in $\Psi$DO formulation}

We now demonstrate the equivalence between the PSF method and the $\Psi$DO formulation (\ref{equ:lr-psido}) using a low-rank symbol representation. Specifically, each PSF $p_k(\vec x)$ can be converted into a row of the symbol by taking its Fourier transform:
\begin{equation}
    s(\vec x_k, \vec \xi) = (2\pi)^2 \overline{\widehat{p}}_k(\vec \xi).
\end{equation}
Too see this, notice that $\mathcal{H}_d$ is self-adjont and hence we have
\begin{align}
    &p_k(\vec x) = \mathcal{H}_d^* \delta_{\vec x_k}(\vec x_k + \vec x) 
    = \overline{\mathcal{H}_d \delta_{\vec x_k + \vec x}(\vec x_k)} \nonumber\\
    &= \int e^{-i\vec x_k \cdot \vec \xi} \overline{s(\vec x_k, \vec \xi)} \overline{\widehat{\delta}_{\vec x_k + \vec x}(\vec \xi)} d\vec \xi \nonumber\\
    &= \frac{1}{(2\pi)^2}\int e^{-i\vec x_k \cdot \vec \xi} \overline{s(\vec x_k, \vec \xi)} e^{i(\vec x_k + \vec x)\cdot\vec \xi} d\vec \xi \nonumber\\
    &= \frac{1}{(2\pi)^2} \int e^{i\vec x \cdot \vec \xi} \overline{s(\vec x_k, \vec \xi)} d\vec \xi \nonumber\\
    &= \frac{1}{(2\pi)^2} \mathcal{F}^{-1} \left[\overline{s(\vec x_k, \vec \xi)} \right](\vec x),
\end{align}
which implies
\begin{equation}
s(\vec x_k, \vec \xi) 
% = (2\pi)^2 \mathcal{F}[p_k(-\vec x)](\vec \xi)
% = (2\pi)^2 \widehat p_k(-\vec \xi) 
= (2\pi)^2 \overline{\widehat{p}}_k(\vec \xi).
\end{equation}

On the other hand, the product-convolution formulation (\ref{equ:psf-pc}) can be rewritten in the $\Psi$DO framework using the convolution theorem and the self-adjointness of $\mathcal{H}_d$:
\begin{equation}
    \mathcal{H}_d 
    \approx (2\pi)^2 \sum_{k=1}^r \mathcal{D}_{w_k} \mathcal{F}^{-1} \mathcal{D}_{\overline{\widehat{p}}_k} \mathcal{F},
\end{equation}
\begin{equation}
    \Leftrightarrow s(\vec x, \vec \xi) 
    \approx (2\pi)^2 \sum_{k=1}^r w_k(\vec x) \overline{\widehat{p}}_k(\vec \xi)
    = \sum_{k=1}^r w_k(\vec x) s(\vec x_k, \vec \xi).
\end{equation}
Therefore, the PSF method can be interpreted as computing rows of the symbol and interpolating them to approximate the full symbol.

\subsection{Novel low-rank symbol method (PSF+)}

\begin{figure}
    \centering
    \includegraphics[width=\linewidth]{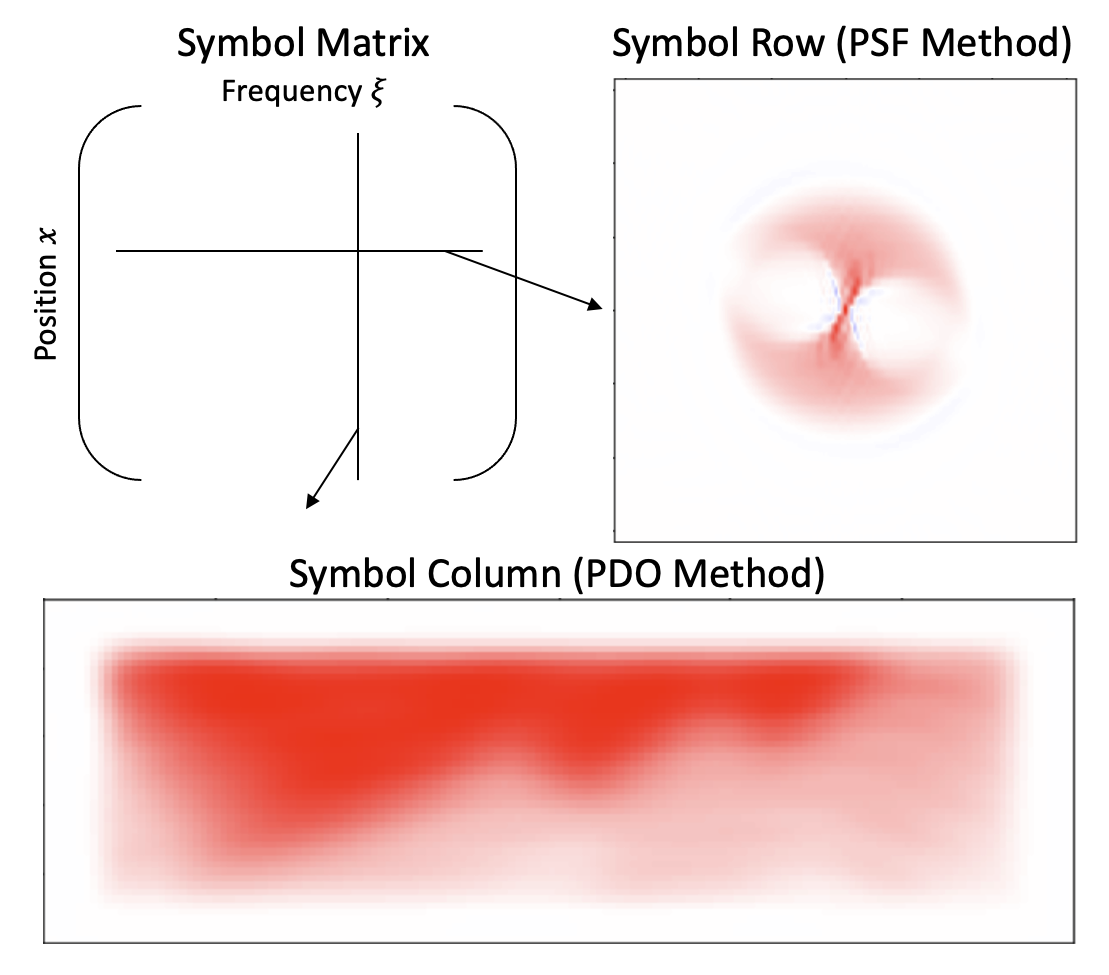}
    \caption{Illustration of the duality between the PSF and PDO methods. The PSF method approximates the misfit Hessian by interpolating sampled rows of the symbol. The PDO method approximates the misfit Hessian by interpolating sampled columns of the symbol.}
    \label{fig:symbol}
\end{figure}

Based on the investigation of the PSF and PDO methods, we found that they share the same underlying intuition, though in a transposed manner. 
While the PSF method leverages locality in spatial space to compute multiple symbol rows through application of the misfit Hessian on delta functions, and subsequently constructs the symbol via interpolation in spatial space, the PDO method operates in a dual manner. It utilizes locality in frequency space to compute multiple symbol columns through application of the misfit Hessian on sinusoids, and then constructs the symbol through interpolation in frequency space. An illustration of this intuition this shown in Fig.~\ref{fig:symbol}.

While both methods achieve efficiency by utilizing locality in spatial or frequency space, they exhibit limitations in practice. In the PSF method, the selected PSFs must be sufficiently spaced to avoid overlap. As a result, the associated weighting functions may fail to accurately capture spatial variations in the PSFs, particularly when their magnitude and shape change rapidly due to sharp contrasts in the model parameters.
On the other hand, the PDO method's interpolation relies on specific conditions---such as smoothness and sufficient frequency coverage---that are often violated in practical applications. 
For instance, due to the lack of high-frequency information in practice, the symbol $s(\vec x, \vec \xi)$ should approach zero for large $|\vec \xi|$, which conflicts with the linear extension assumption.
To address these limitations, we propose a low-rank symbol approximation that simultaneously manipulates both symbol rows and columns. This approach leverages the PSF method's accuracy in spatial space to enhance the PDO method's limitations in frequency interpolation, and vice versa, thereby improving overall symbol estimation.

Suppose the approximated symbol rows $\widetilde{S}[I_r, :]$ and symbol columns $\widetilde{S}[:, I_c]$ are precomputed and saved, where $I_r$ and $I_c$ are the indices of sampled rows and columns, respectively. Here, $\widetilde{S}[I_r, :]$ denotes the submatrix of $\widetilde{S}$ consisting of the rows indexed by $I_r$ and all columns, and $\widetilde{S}[:, I_c]$ denotes the submatrix consisting of all rows and the columns indexed by $I_c$. We design the PSF+ method using the following steps:
\begin{enumerate}
    \item Start with the same strategy as the PSF method by computing interpolation weights. This results in a low-rank approximation of the symbol, $S_{\mathrm{psf}} = W \widetilde{S}[I_r, :]$, where $W$ is the weighting matrix. 
    \item Apply the singular value decomposition (SVD) to the computed symbol rows $\widetilde{S}[I_r, :]$ to extract the dominant components (this technique to extract information from PSFs via SVD has a long history, for example, in \citet{denis2015fast}). This step forms a basis matrix $B$ for the row space of the symbol, yielding an updated expression, $S_{\mathrm{psf}} = A_0 B$, with reduced rank, effectively eliminating redundant information from the PSFs.
    \item Refine the matrix $A_0$ by minimizing the difference between the updated symbol columns and the computed symbol columns. Tikhonov regularization is applied to prevent significant deviations and maintain consistency among the symbol rows. The resulting small-sized least-squares problem is
    \begin{equation}\label{equ:psfplus}
        \min_{A\in\mathbb{R}^{N\times r}}  \|AB[:,I_c] - \widetilde{S}[:, I_c]\|_F^2 + \alpha \|A-A_0\|_F^2,
    \end{equation}
    where the regularization coefficient $\alpha = \|B[:,I_c]\|_F^2 / \|B\|_F^2$ balances the refinement term and the penalty term.
\end{enumerate}

This novel approach leverages the strengths of both the PSF and PDO methods, resulting in a more accurate low-rank symbol approximation that can be converted into a fast $\Psi$DO approximation of $H_d$. 

\subsection{Practical corrections for the limitations of the methods}

Both the PSF and PDO methods are designed to operate on the high frequency part of $H_d$. In the PSF method, low-frequency components compromise the locality of point spread functions, while the PDO method requires sampling sufficiently high frequencies for spectral separation. To handle this, we apply a high-pass filter to $H_d$ before the approximation, removing low-frequency components that are beyond the scope of these approaches. Consequently, the proposed methods approximate the symbol $s(\vec x, \vec \xi)$ only for high-frequency part, leaving the low-frequency portion unaddressed.

Fortunately, the low-frequency part of the $H_d$ is smooth and hence low-rank. This enables efficient correction methods using the low-rank approach described in Section 3.1, which requires us to apply $H_d$ to only a few vectors. 
In deterministic seismic inversion, low-rank correction is typically unnecessary. 
If the Hessian approximation (i.e., the matrix $\widetilde{H} := \widetilde{H}_d + R \approx H$) is used as an initial guess for the Hessian in a quasi-Newton method such as L-BFGS, the low-rank Hessian updates in the quasi-Newton method automatically correct errors in the low frequency part of the Hessian approximation. If the Hessian approximation is used as a preconditioner in a Newton-Krylov method, the Krylov iteration automatically corrects for errors in the low frequency part of the Hessian approximation.
In contrast, for UQ, precise Hessian estimation is essential for sampling from the posterior distribution. Although there are a small number of eigenvalues corresponding to the low frequency part of the Hessian, these eigenvalues are large and therefore significantly influence samples. Thus, low-rank correction is required to refine the Hessian approximation before using it for posterior sampling.

For UQ, the Hessian approximation must be symmetric positive definite, as it represents the inverse of the covariance matrix in a Gaussian distribution. Therefore, rather than directly approximating the symbol, an approach that may yield a non-symmetric operator, we instead approximate the square root of $H_d$. The overall approximation procedure remains unchanged, except that we apply a pointwise square root to the computed symbol rows or columns, discarding any negative values. This approach is justified by the property that the symbol of the product of two $\Psi$DOs is approximately equal to the pointwise product of their individual symbols. This idea is not new; it is also employed in the PSF deconvolution method \citep{yang2022efficient} , where the pointwise inverse of symbol rows is used to construct an approximate inverse Hessian.

\section{Numerical Results}

We implement the proposed methods to accelerate UQ for an ideal quadratic model and a modified Marmousi Model. Throughout the numerical results, LR, PSF, PDO, and PSF+ denote the low-rank, point spread function, pseudo-differential operator probing, and novel low-rank symbol Hessian approximation methods, respectively.

\subsection{Ideal quadratic model}

\begin{figure}
    \centering
    \begin{subfigure}{\linewidth}
        \centering
        \includegraphics[width=\textwidth]{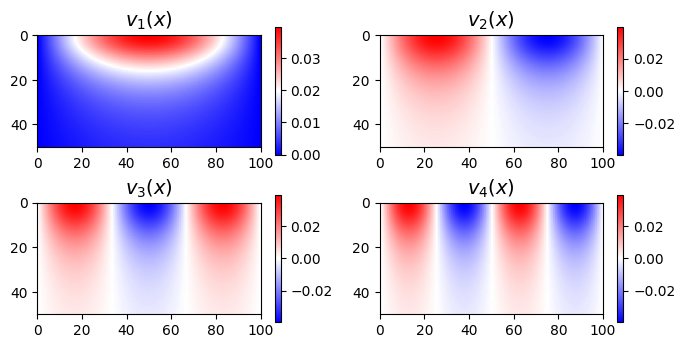}
        \caption{The leading modes $v_k(x)$ (columns of $V_r$) of the low-rank part of the Hessian.}\label{fig:toy-hess-lr}
    \end{subfigure}
    \begin{subfigure}{\linewidth}
        \centering
        \includegraphics[width=\textwidth]{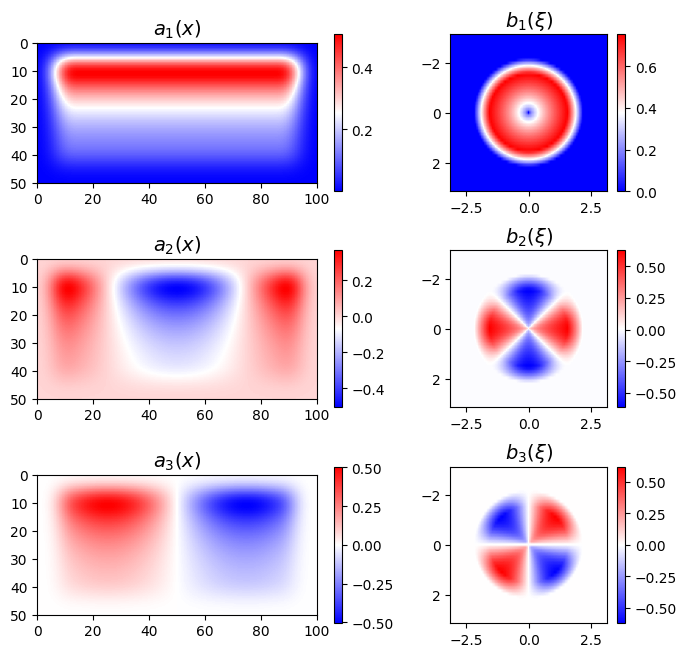}
        \caption{The components $a_k$ and $b_k$ of the $\Psi$DO $A$, which forms the $\Psi$DO part of the idealized Hessian.}\label{fig:toy-hess-psido}
    \end{subfigure}
    \caption{The components of the Hessian in the ideal quadratic model.}
    \label{fig:toy-hess}
\end{figure}

Since seismic inversion is large-scale, ill-posed and highly nonlinear, the true uncertainty is unavailable. We  construct an ideal quadratic model where the exact uncertainty is available, allowing a straightforward evaluation of our method's performance in UQ. This quadratic model problem mimics the seismic inversion problem. The misfit Hessian of the model problem,
\begin{equation}\label{equ:toy-hess}
    H_d = AA^T + V_r \Lambda_r V_r^T,
\end{equation}
is explicitly constructed as an order-1 $\Psi$DO $AA^T$, which mimics the reflection part of the misfit Hessian in seismic inversion, plus a low-rank component $V_r \Lambda_r V_r^T$, which mimics the low-frequency part of the misfit Hessian in seismic inversion. 
Specifically, $\Lambda_r$ is a $10 \times 10$ diagonal matrix with positive diagonals, and $V_r$ has 10 orthonormal columns, each representing a plane wave with high amplitude near the surface that decays rapidly with depth (see Fig.~\ref{fig:toy-hess-lr}). 
The $\Psi$DO $A$ is of order-0.5 (so that $AA^T$ is of order-1) and has a low-rank symbol that diminishes with depth:
\begin{equation}\label{equ:toy-symbol}
    s_A(\vec x, \vec \xi) = w(\vec x) P_{0.5}(|\vec\xi|) \frac{0.1 + \cos^2(\arg \vec \xi + \pi x/x_{\mathrm{max}})}{0.1+(z/z_{\mathrm{max}})^2}.
\end{equation}
The window function $w(\vec x)$ suppresses values near the boundaries.
The function $P_{0.5}(|\vec\xi|)$ asymptotically behaves as $|\vec\xi|^{1/2}$ but decays to zero for large $|\vec\xi|$, reflecting realistic loss of out-of-band frequency information. 
The denominator scales the symbol by depth $z$, while the cosine term accounts for variations in reflector orientation. 
This symbol follows the low-rank format (\ref{equ:lr_symbol}), with a rank-3 decomposition using the compound angle formula:
\begin{equation}
    s_A(\vec x, \vec \xi) = \sum_{k=1}^3 a_k(\vec x) b_k(\vec \xi),
\end{equation}
where
\begin{align}
&a_1(\vec x) = \frac{w(\vec x)}{0.1+(z/z_{\mathrm{max}})^2}, \quad b_1(\vec \xi) = 0.6 \times P_{0.5}(|\vec\xi|),\\
&a_2(\vec x) = \frac{w(\vec x) \cos(2\pi x/x_{\mathrm{max}})}{0.1+(z/z_{\mathrm{max}})^2},\\
&b_2(\vec \xi) =  \frac12 \cos(2\arg\vec\xi) P_{0.5}(|\vec\xi|),\\
&a_3(\vec x) = \frac{w(\vec x) \sin(2\pi x/x_{\mathrm{max}})}{0.1+(z/z_{\mathrm{max}})^2},\\
&b_3(\vec \xi) = -\frac12 \sin(2\arg\vec\xi) P_{0.5}(|\vec\xi|),
\end{align}
% \[
% s_A(\vec x, \vec \xi) = \frac{w(\vec x)}{0.1+(z/z_{\mathrm{max}})^2} \times 0.6 P_{0.5}(|\vec\xi|) 
% \]
% \[
% + \frac{w(\vec x) \cos(2\pi x/x_{\mathrm{max}})}{0.1+(z/z_{\mathrm{max}})^2} \times \frac12 \cos(2\arg\vec\xi) P_{0.5}(|\vec\xi|) 
% \]
% \begin{equation}
%     - \frac{w(\vec x) \sin(2\pi x/x_{\mathrm{max}})}{0.1+(z/z_{\mathrm{max}})^2} \times \frac12 \sin(2\arg\vec\xi) P_{0.5}(|\vec\xi|),
% \end{equation}
where each term is shown in Fig.~\ref{fig:toy-hess-psido}.

\begin{figure*}
    \centering
    \includegraphics[width=0.49\textwidth]{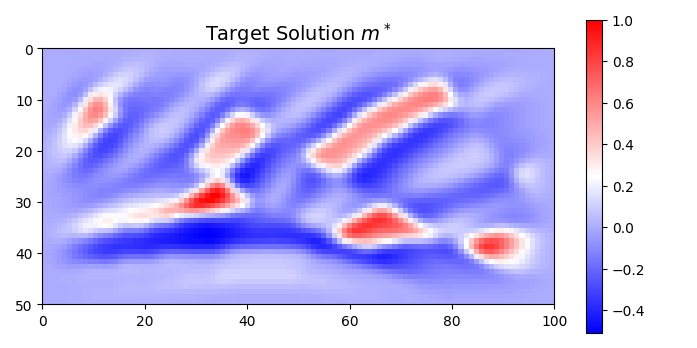}
    \includegraphics[width=0.49\textwidth]{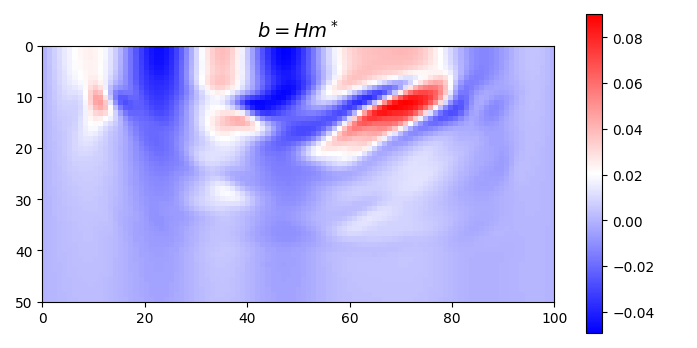}
    \caption{The target solution and the right-hand-side in the ideal quadratic model.}
    \label{fig:toy-model-xb}
\end{figure*}

The misfit is defined in a quadratic form:
\begin{equation} \label{equ:toy-model}
    \Phi_d(\vec m) = \frac{1}{2} \vec m^T H_d \vec mx - \vec b^T \vec m + \frac{1}{2} \vec b^T \vec m^*,
\end{equation}
where $\vec b = H_d \vec m^*$, with target solution $\vec m^*$ constructed as a normalized difference between a subset of the Marmousi model and its smoothed variant (see Fig.~\ref{fig:toy-model-xb}). 
This misfit function resembles seismic misfit and reaches a global minimum of zero at $\vec m^*$.

% \begin{figure*}
%     \centering
%     \includegraphics[width=\textwidth]{images/toy-model.png}
%     \caption{Ideal quadratic model settings. Left: The leading modes in the low-rank part. Middle: The functions $a_k$ of the ideal misfit Hessian's $\Psi$DO part. Right: The functions $b_k$ of the ideal misfit Hessian's $\Psi$DO part.}
%     \label{fig:toy-model}
% \end{figure*}

To mitigate boundary artifacts, we apply a window function $W$ that suppresses values near the domain edges, which is particularly suited to seismic applications. Near-surface structures are typically well-recovered by inexpensive methods, making inversion unnecessary in this region. Furthermore, the left, right, and bottom boundaries are poorly informed, yielding unreliable results in these areas.
Finally, we apply a biharmonic prior (see Section 2), resulting in the objective function and posterior distribution defined as
\begin{equation}
    \Phi(\vec m) = \Phi_d(W \vec m) + \frac{1}{2} \|\vec m\|^2_{\prcov^{-1}},
\end{equation}
and 
\begin{equation}
    \pi_{\mathrm{post}}(\vec m) \propto \exp \left\{ -\Phi_d(W \vec m) - \frac{1}{2} \|\vec m\|^2_{\prcov^{-1}} \right\}.
\end{equation}

\subsubsection{Hessian approximation}

We approximate the Hessian using the methods discussed. 

\begin{figure*}
    \centering
    \includegraphics[width=0.5\textwidth]{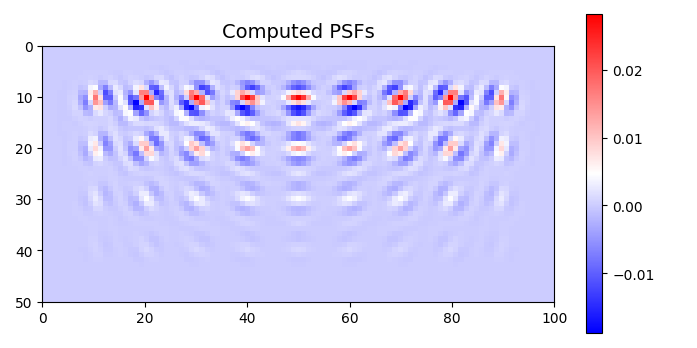}
    \caption{Computed PSFs in the ideal quadratic model.}
    \label{fig:toy-app-psf-vec}
\end{figure*}

In the PSF method, we compute 36 PSFs by applying $H_d$ to 4 vectors. 
Fig.~\ref{fig:toy-app-psf-vec} displays the computed PSFs, highlighting their locality properties. Moreover, the orientation of each PSF varies across the domain---an effect governed by the cosine term in the symbol formulation~(\ref{equ:toy-symbol}) of the $\Psi$DO $A$. These PSFs are converted to symbol rows, and the weighting functions are created accordingly. 

\begin{figure*}
    \centering
    \includegraphics[width=0.65\textwidth]{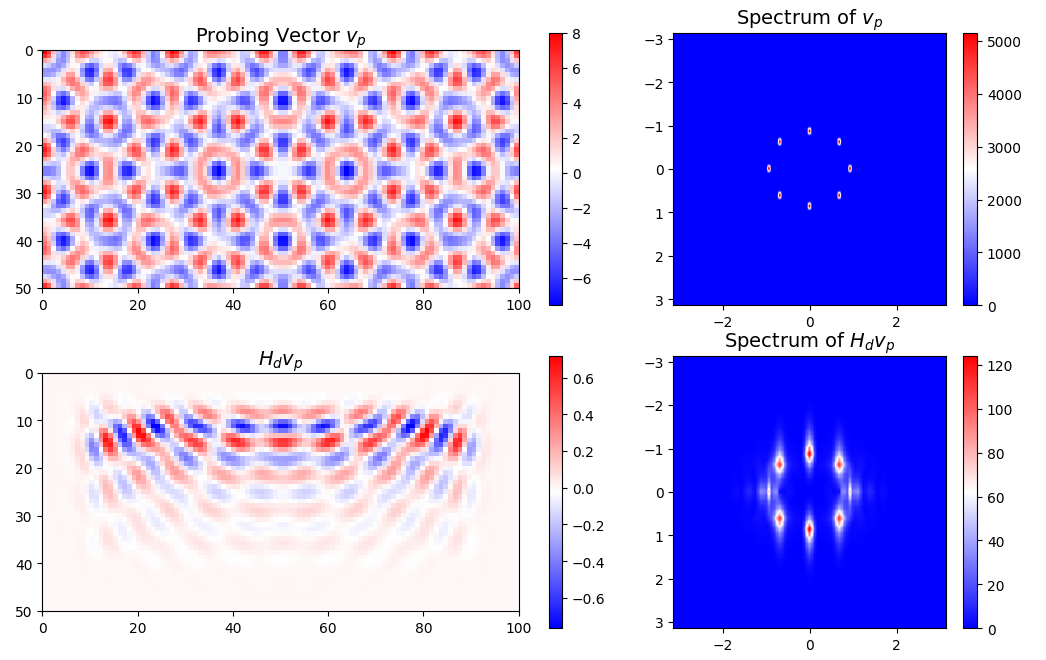}
    \caption{The PDO method in the ideal quadratic model: The probing vector and the result of applying misfit Hessian to the probing vector.}
    \label{fig:toy-app-pdo-vec}
\end{figure*}

In the PDO method, we compute 8 symbol columns by applying $H_d$ to one probing vector. Fig.~\ref{fig:toy-app-pdo-vec} shows the probing vector, the result of applying $H_d$ to the probing vector, and their Fourier transforms, illustrating the locality of $H_d$ in frequency space. The weighting functions combine Fourier interpolation over the angle and linear extension over the radius. 
% In addition, we filtered out the out-of-band high-frequency part of the weighing function, in order to correct the approximation in high-frequencies. 

The PSF+ method utilizes both symbol rows and columns computed from the PSF and PDO method.
The row basis $b_k(\vec \xi)$ is first extracted from the symbol rows via SVD, where the leading three singular vectors sufficiently span the row space of the symbol.
Then, the functions $a_k(\vec x)$ are computed by solving the small-size optimization problem (\ref{equ:psfplus}).

\subsubsection{Deterministic inversion}

To evaluate the performance of the preconditioners derived from Hessian approximations, we compare L-BFGS convergence with two baselines: the original L-BFGS algorithm and L-BFGS with a prior-preconditioner. 
In addition to traditional metrics such as data misfit and gradient norm of the objective, we plot the $L^2$-norm of the solution error to assess the overall quality of the computed model relative to the ground truth, where the updates in deeper regions are less effectively captured by the misfit or gradient.

\begin{figure*}
    \centering
    \includegraphics[width=0.33\textwidth]{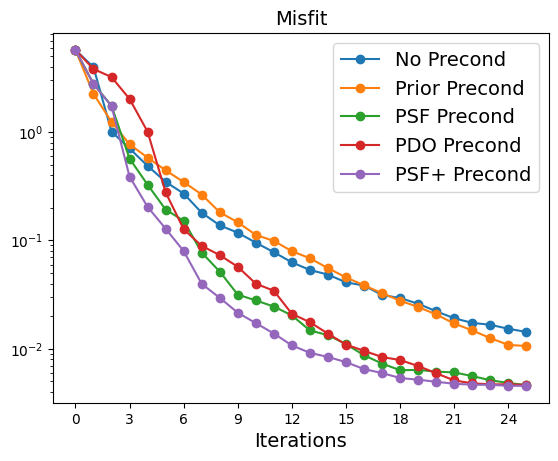}%
    \includegraphics[width=0.33\textwidth]{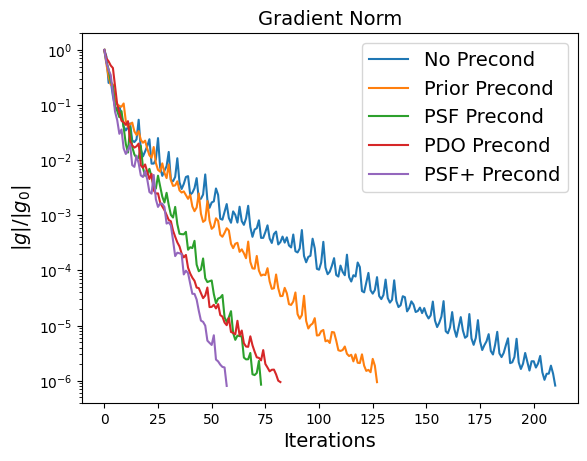}%
    \includegraphics[width=0.33\textwidth]{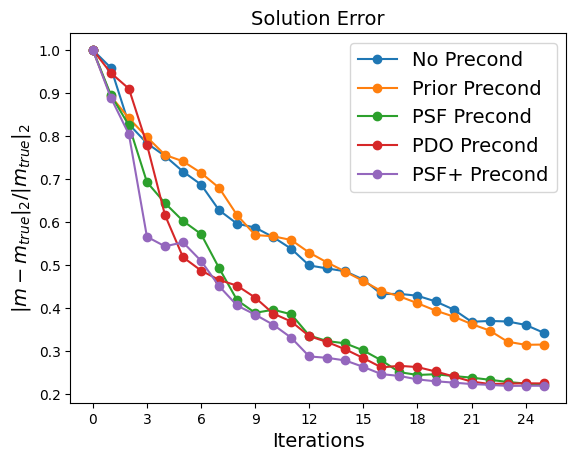}%
    \caption{L-BFGS convergence in the ideal quadratic model.}
    \label{fig:toy-map}
\end{figure*}

\begin{figure*}
    \centering
    \includegraphics[width=\textwidth]{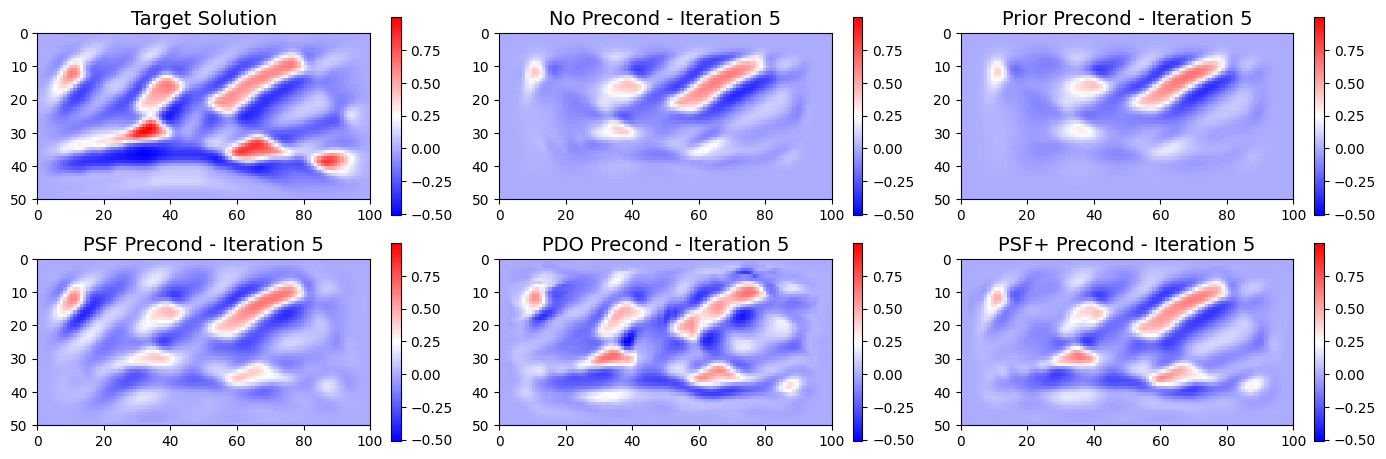}
    \caption{L-BFGS solutions at the 5-th iteration in the ideal quadratic model.}
    \label{fig:toy-map-sol}
\end{figure*}

% From a short view
Examining the misfit curve in Fig.~\ref{fig:toy-map}, we observe that the three Hessian-based preconditioners lead to slower initial misfit decay in the first 5 iterations due to missing low-frequency components. 
However, this low-frequency part is naturally addressed by the L-BFGS algorithm, leading to faster misfit decay in later iterations for Hessian-based methods. 
A similar trend is observed in the solution error plot in Fig.~\ref{fig:toy-map}, where the benefits of Hessian-based preconditioners are more evident due to better deep part recovery. 
Additionally, Fig.~\ref{fig:toy-map-sol} shows the solutions at the fifth iteration, where the high-rank Hessian approximations notably improve deep-region recovery compared to baselines.

% From a long view
For long-term convergence, the gradient norms in Fig.~\ref{fig:toy-map} show that using the Hessian preconditioners with L-BFGS yields a faster convergence rate. This is further confirmed in \autoref{tab:toy-map}, which lists the iteration counts required to reach convergence, defined as a $10^{-6}$ reduction in gradient norm.

\begin{table}
    \centering
    \caption{L-BFGS Computational Cost in Ideal Quadratic Model}
    \label{tab:toy-map}
    \begin{tabular}{lc}
        \toprule
                             & L-BFGS Steps \\ \hline
        No Preconditioner    &     210      \\
        Prior Preconditioner &     127      \\
        PSF                  &      73\\
        PDO                  &      82\\
        PSF+                 &      57\\ \hline
    \end{tabular}
\end{table}

\subsubsection{Uncertainty quantification}

In addition to the Hessian approximations, we apply a rank-100 correction using randomized SVD. This enhanced Hessian approximation is then used to construct a proposal distribution for sampling from the Laplace approximation of the posterior. We generate 50,000 samples using the MCMC-gpCN algorithm. For comparison, two baselines are constructed: one using the MCMC-pCN algorithm and the other using the MCMC-gpCN with only the low-rank Hessian approximation with 100 leading eigenvalues and eigenvectors. 
% This setup highlights the benefits of our high-rank Hessian approximations.

\begin{figure*}
    \centering
    \includegraphics[width=0.5\textwidth]{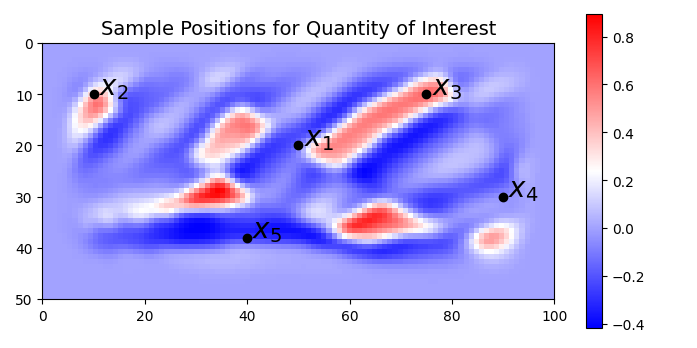}
    \caption{Sample positions for MCMC evaluations in the ideal quadratic model.}
    \label{fig:toy-uq-qoi}
\end{figure*}

\begin{figure*}
    \centering
    \includegraphics[width=\textwidth]{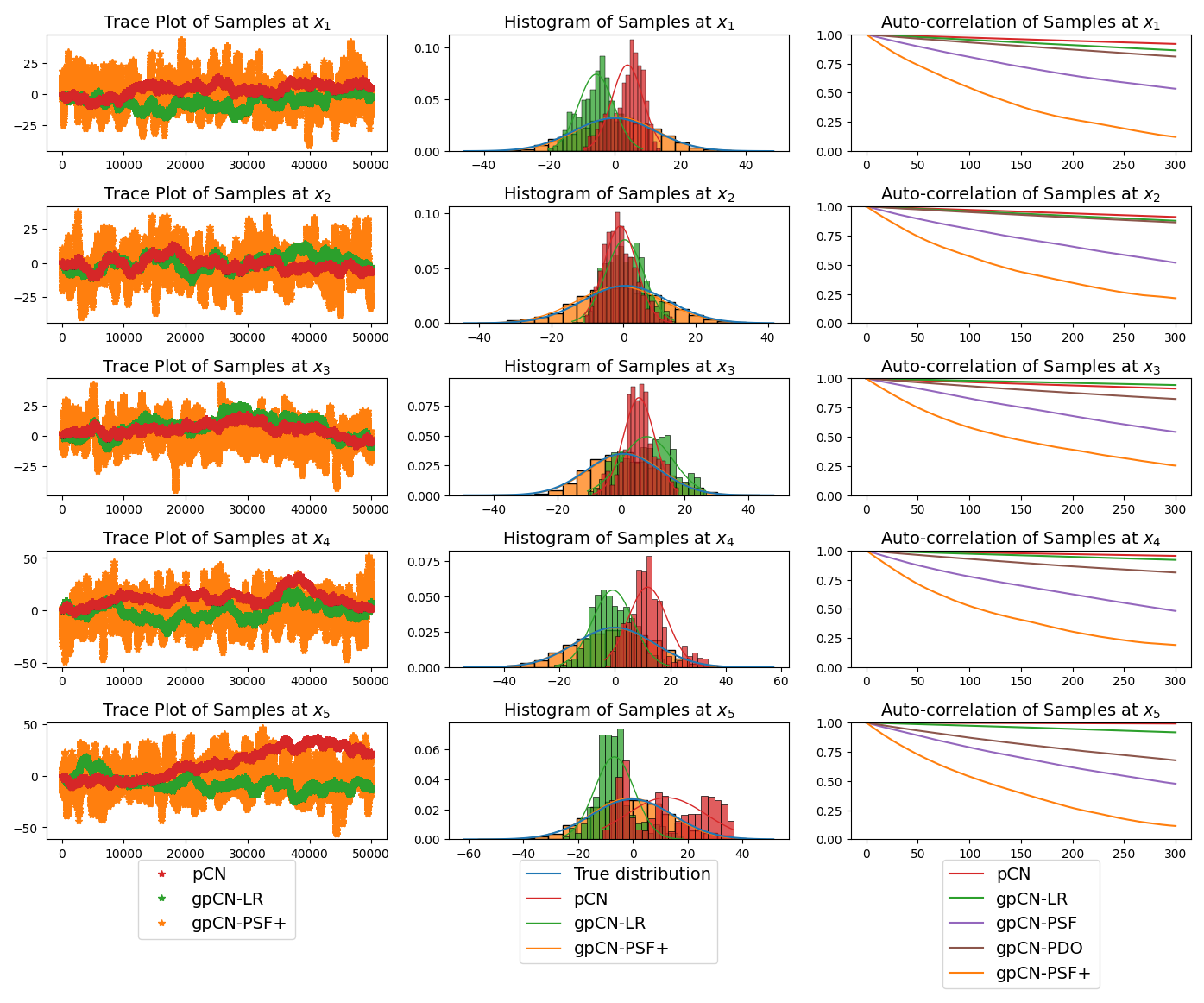}
    \caption{Performance of MCMC in the ideal quadratic model.}
    \label{fig:toy-uq}
\end{figure*}

To demonstrate improvements in UQ, we examine sample values at five fixed points (Fig.~\ref{fig:toy-uq-qoi}). 
MCMC performance is evaluated with three metrics in Fig.~\ref{fig:toy-uq}. 
First, the sample traces show that our method achieves superior mixing compared to traditional approaches. 
Second, histograms reveal that our method captures accurate distributions, while both pCN and gpCN-LR restrict samples to narrower ranges, underestimating variation. 
Third, the autocorrelation plot shows faster decay with our PSF+ method, indicating more efficient sampling. 
\autoref{tab:toy-uq} further reports effective sample sizes across sample chains at different points, confirming the PSF+ method’s performance advantage.

\begin{table}
    % \centering
    \caption{MCMC for Ideal Quadratic Model}
    \label{tab:toy-uq}
    \begin{tabular}{lccccc}
        \hline
                   &           \multicolumn{5}{c}{Effective Sample Size}            \\
                   & $\vec x_1$ & $\vec x_2$ & $\vec x_3$ & $\vec x_4$ & $\vec x_5$ \\ \hline
        pCN        &    7.02    &   11.50    &    6.01    &   10.08    &    3.00    \\
        gpCN-LR         &    7.89    &   10.55    &    6.93    &    9.99    &    7.65    \\
        gpCN-PSF     &   57.48    &   55.53    &   23.86    &   51.43    &   45.32    \\
        gpCN-PDO     &   24.20    &   14.00    &   16.93    &   24.36    &   42.33    \\
        gpCN-PSF+ &   132.66   &   104.42   &   59.82    &   117.95   &   109.23   \\ \hline
    \end{tabular}
\end{table}

\begin{figure*}
    \centering
    \includegraphics[width=\textwidth]{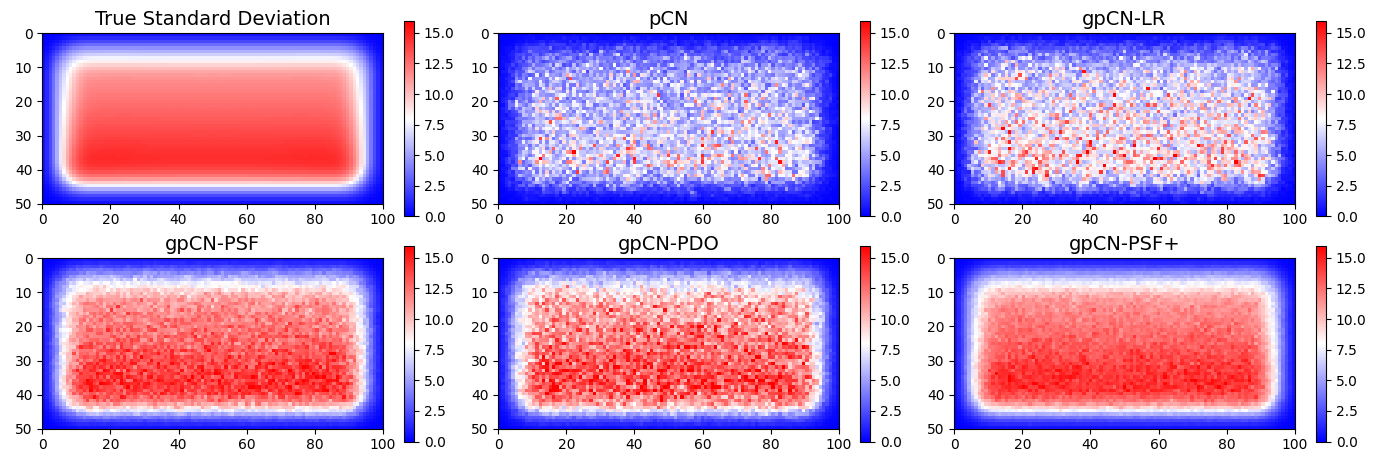}
    \caption{Standard deviations computed from the MCMC samples with different Hessian approximations in the ideal quadratic model.}
    \label{fig:toy-uq-std}
\end{figure*}

Finally, Fig.~\ref{fig:toy-uq-std} compares standard deviations from MCMC samples to the true standard deviation from the inverse Hessian. In the pCN and gpCN-LR baselines, standard deviations are generally underestimated, with noisy results suggesting non-convergent sampling and unreliable statistics. In contrast, our Hessian approximation captures higher-frequency modes, yielding accurate standard deviations and smoother results, demonstrating improved MCMC convergence. Among the methods, the PSF+ approach provides the most accurate and smoothest standard deviations, closely matching the ground truth.

\subsection{Marmousi model}

The Marmousi model \citep{brougois1990marmousi} is a widely used benchmark in seismic imaging, featuring a complex 2D geological structure with significant velocity variations both horizontally and vertically.

\begin{figure*}
    \centering
    \includegraphics[width=0.49\textwidth]{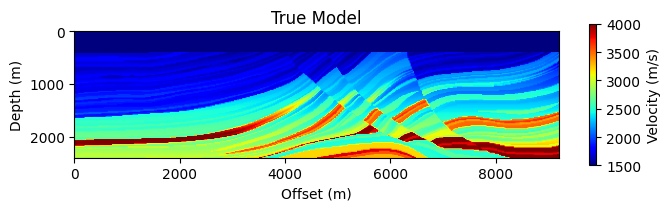}
    \includegraphics[width=0.49\textwidth]{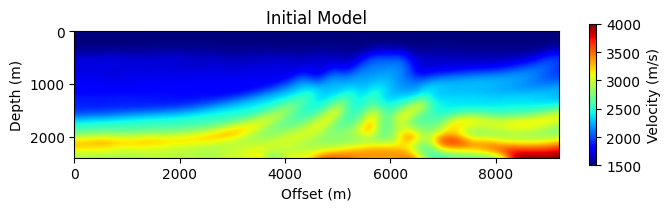}
    \caption{Ground truth and the initial guess in the Marmousi model.}
    \label{fig:mar-model}
\end{figure*}

We made several modifications to the model: trimming the bottom 1000m to remove a highly ill-conditioned region and adding a 400m-deep water layer at the top. The resulting model spans 9200m$\times$2400m and is discretized on a 461$\times$121 grid. We place 92 sources at 100m intervals along the top boundary and 153 receivers at 60m intervals. Fig.~\ref{fig:mar-model} shows both the ground truth velocity model and an initial model, which is a smoothed version of the ground truth. Wave simulations are performed using a finite difference method implemented in the programming language Julia \citep{bezanson2017julia}.
%An region-extension method is applied to address the reflections on boundaries.

\subsubsection{Deterministic inversion}

\begin{figure*}
    \centering
    \includegraphics[width=0.5\textwidth]{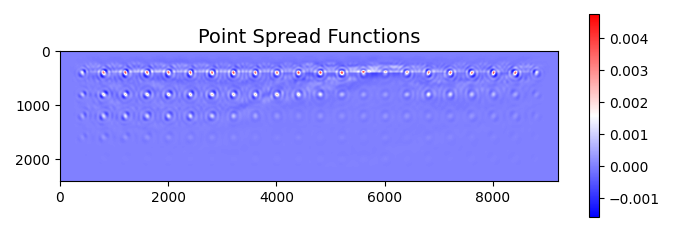}
    \caption{Computed PSFs for the initial misfit Hessian in the Marmousi model.}
    \label{fig:mar-app-psf-vec}
\end{figure*}

\begin{figure*}
    \centering
    \includegraphics[width=0.8\textwidth]{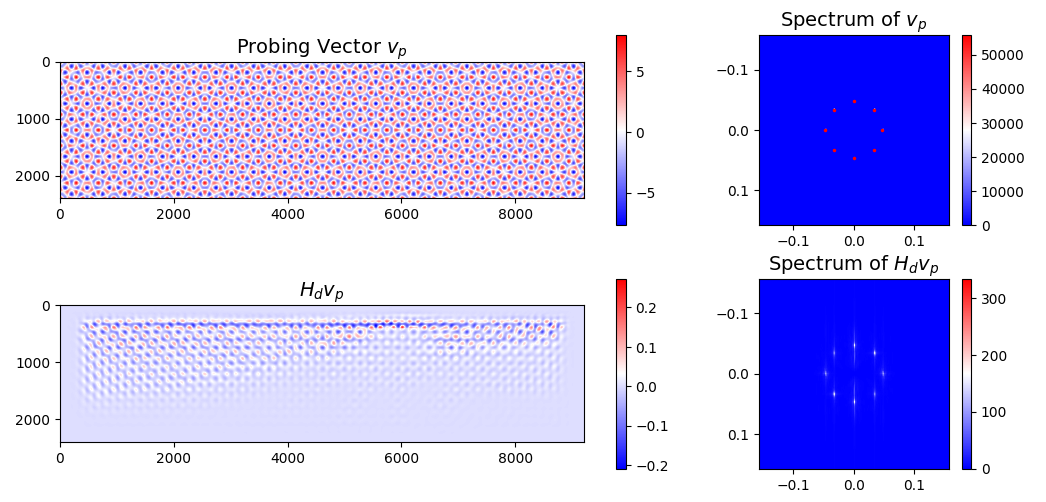}
    \caption{PDO Method: The probing vector and the application of the misfit Hessian to the probing vector.}
    \label{fig:mar-app-pdo-vec}
\end{figure*}

The Hessian is approximated at the starting model using the three proposed high-rank approximation methods. Computed PSFs are shown in Fig.~\ref{fig:mar-app-psf-vec}. The PDO probing vector and the result of applying the misfit Hessian to the probing vector are shown in Fig.~\ref{fig:mar-app-pdo-vec}. These figures validate the locality property of the seismic misfit Hessian in both spatial and frequency spaces.

\begin{figure*}
    \centering
    \includegraphics[width=\textwidth]{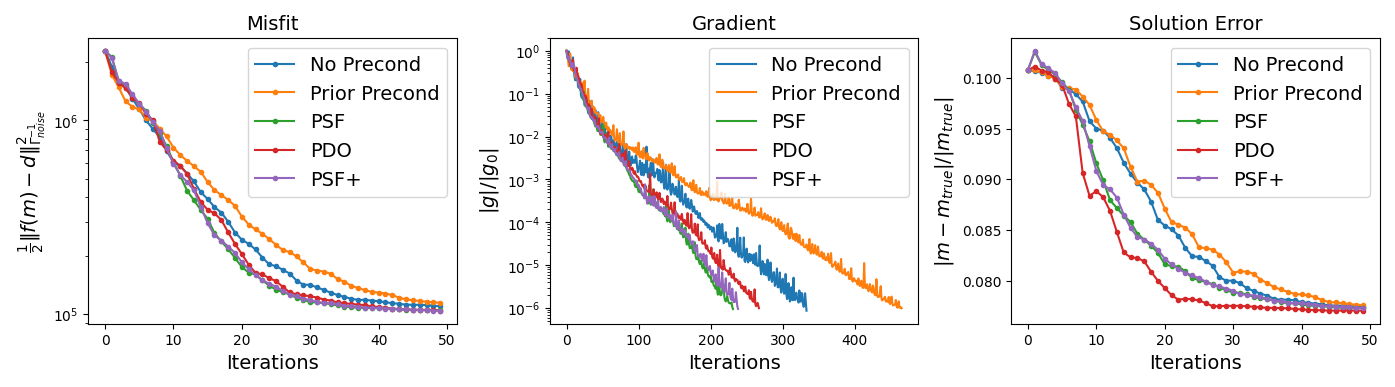}
    \caption{Comparison of preconditioners in L-BFGS solvers for the Marmousi model. The Hessian-based preconditioners lead to faster decay of data misfit, gradient norm, and solution error.}
    \label{fig:mar-map}
\end{figure*}

\begin{figure*}
    \centering
    \includegraphics[width=\textwidth]{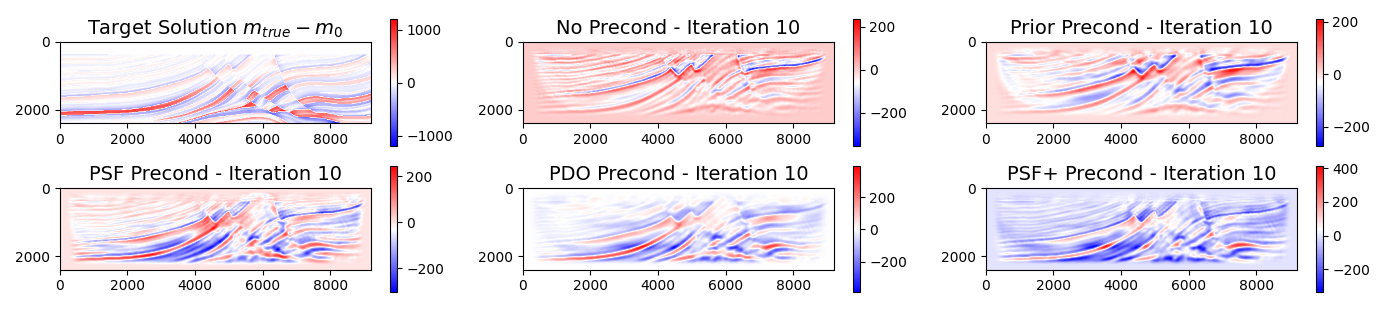}
    \caption{Comparison of preconditioners in L-BFGS solvers for the Marmousi model. Each picture shows the difference between the current velocity model and the initial model. The Hessian-based preconditioners recover details faster.}
    \label{fig:mar-map-sol}
\end{figure*}

Examining the misfit curve in Fig.~\ref{fig:mar-map}, the Hessian-based preconditioners achieve faster misfit reduction compared to the baselines, though the improvement is modest. In contrast, the solution error curves show more pronounced benefits. Fig.~\ref{fig:mar-map-sol} displays updated models at the tenth iteration, highlighting enhanced model recovery in terms of both structure and magnitude. \autoref{tab:mar-map} confirms that the Hessian-based preconditioners significantly accelerate convergence and reduce computation time.

\begin{table}
    \centering
    \caption{L-BFGS Computational Cost for Marmousi Model}
    \label{tab:mar-map}
    \begin{tabular}{lc}
        \hline
                                & L-BFGS Steps  \\ \hline
        No Preconditioner       & 333           \\
        Prior Preconditioner    & 465           \\
        PSF                     & 231           \\
        PDO                     & 267           \\
        PSF+                 & 238           \\ \hline
    \end{tabular}
\end{table}

\subsubsection{Uncertainty quantification}

The Hessian is re-approximated at the MAP model and corrected by the low-rank method, with 100 leading eigenvalues and eigenvectors, for an accurate Laplace approximation to the posterior.
Unlike the ideal quadratic model, where ground truth uncertainty is precisely known, the Marmousi model lacks an exact uncertainty reference. To address this, we establish a benchmark using a brute-force low-rank approximation of the Hessian at the MAP model, followed by MCMC-gpCN sampling. While computationally impractical for real-world applications, this benchmark provides a valuable reference for evaluating our proposed methods.

\begin{figure*}
    \centering
    \includegraphics[width=0.6\textwidth]{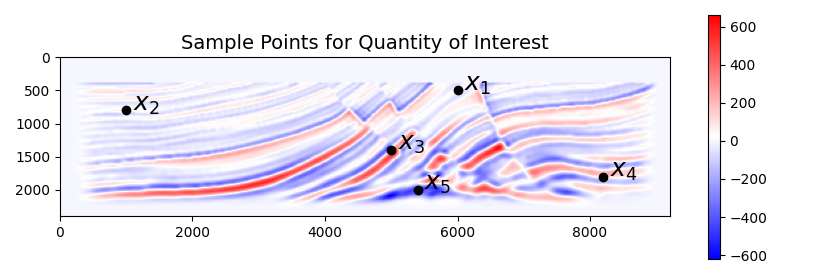}
    \caption{Sample positions for MCMC evaluations in the Marmousi model.}
    \label{fig:mar-uq-qoi}
\end{figure*}

\begin{figure*}
    \centering
    \includegraphics[width=\textwidth]{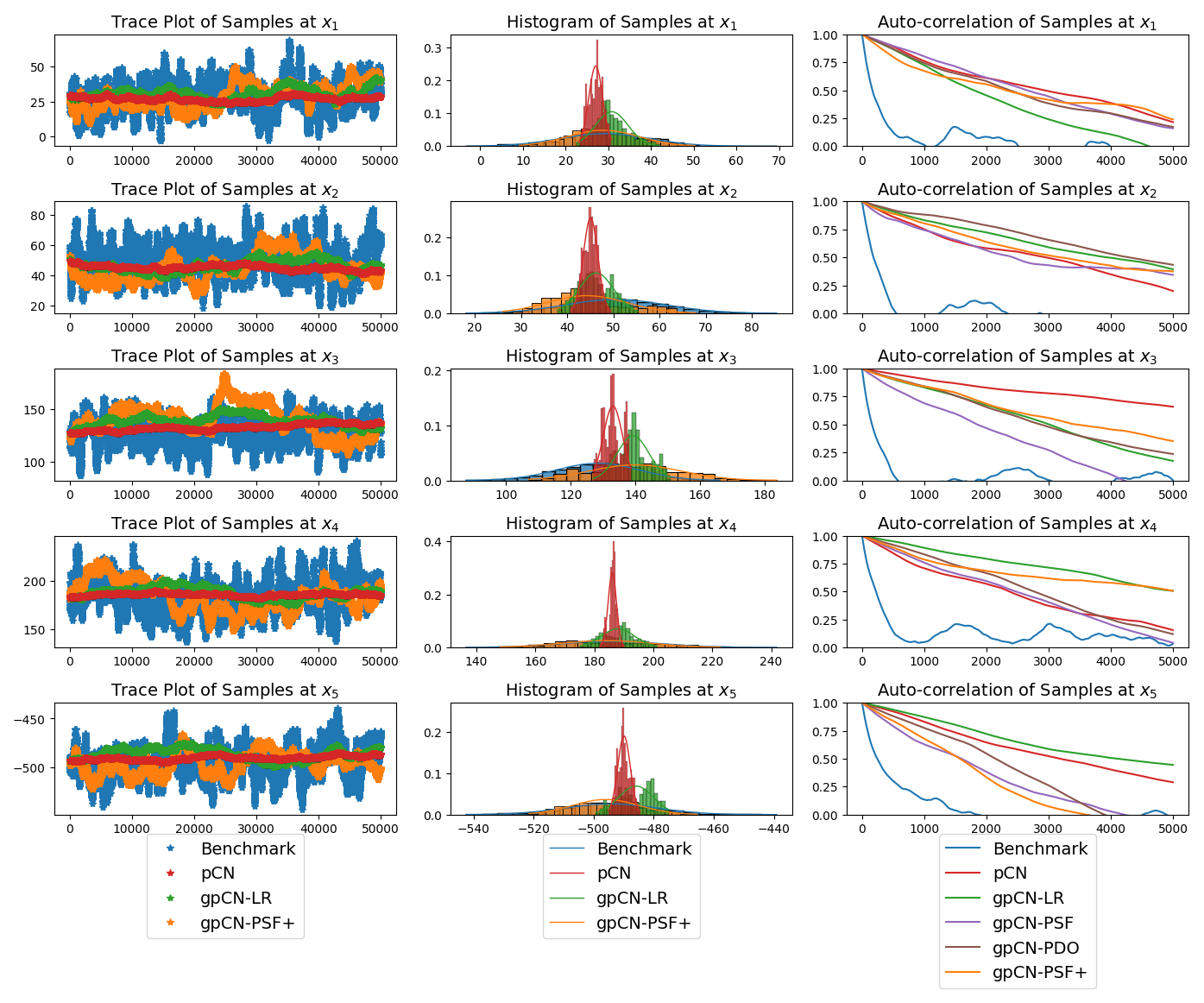}
    \caption{Performance of MCMC in the Marmousi model.}
    \label{fig:mar-uq}
\end{figure*}

MCMC results are evaluated at five fixed points (Fig.~\ref{fig:mar-uq-qoi}) using three metrics, as shown in Fig.~\ref{fig:mar-uq}. Sample traces indicate that the proposed Hessian approximation achieves better sample mixing than conventional methods. Histograms show that the sample distribution from our method closely matches the benchmark, while the pCN and gpCN-LR methods yield narrow, inaccurate distributions, underestimating variation. 

The autocorrelation plot in Fig.~\ref{fig:mar-uq} and the effective sample size (ESS) results reported in \autoref{tab:mar-uq} seem to suggest that the conventional methods (pCN and gpCN-LR) are performing well, when they are actually performing extremely poorly. This is evident, for example, by comparing the histogram of samples at point $x_1$ with the autocorrelation and ESS there. The autocorrelation plot and ESS suggest that gpCN and gpCN-LR are performing roughly as well or better than gpCN-PSF+. However, when we look at the histogram of samples we see that the distributions generated by gpCN and gpCN-LR are highly concentrated and do not match the benchmark distribution, while the distribution generated by gpCN-PSF+ closely matches the benchmark distribution. This mismatch is because the autocorrelation and ESS metrics are heuristic, and do not accurately reflect sampling effectiveness if the proposal distribution is not reflective of the posterior and the MCMC chain has not converged. Although autocorrelation and ESS metrics are common in the literature, they must be used with caution because they can provide false confidence in sampling methods that are performing poorly.

\begin{figure*}
    \centering
    \includegraphics[width=\textwidth]{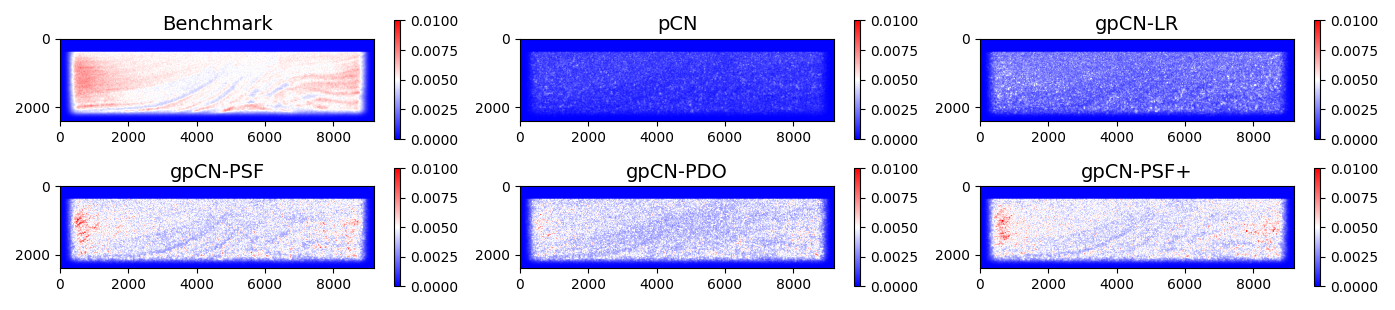}
    \caption{Relative standard deviations computed from the MCMC samples with different Hessian approximations in the Marmousi model.}
    \label{fig:mar-uq-std}
\end{figure*}

Finally, Fig.~\ref{fig:mar-uq-std} presents the standard deviations computed from samples drawn by different methods, scaled by the MAP model velocities. While pCN and gpCN-LR methods underestimate variance with noisy results, our method provides more accurate estimates, demonstrating its effectiveness in capturing true uncertainty in model parameters compared to conventional methods.

\begin{table}
    % \centering
    \caption{MCMC for Marmousi Model}
    \label{tab:mar-uq}
    \begin{tabular}{lccccc}
        \hline
                      &           \multicolumn{5}{c}{Effective Sample Size}            \\
                      & $\vec x_1$ & $\vec x_2$ & $\vec x_3$ & $\vec x_4$ & $\vec x_5$ \\ \hline
    %		gpCN-LR5000     &   69.53    &   100.59   &   110.73   &   27.08    &   53.12    \\
        pCN         &    8.30    &    8.50    &    3.80    &    9.45    &    7.33    \\
        gpCN-LR          &   12.57    &    5.35    &    8.27    &    4.79    &    5.74    \\
        gpCN-PSF      &    6.72    &    6.18    &   12.85    &   10.10    &   12.32    \\
        gpCN-PDO      &    9.23    &    5.04    &    7.84    &    8.94    &   12.10    \\
        gpCN-PSF+  &    7.05    &    6.14    &    6.57    &    4.88    &   16.00    \\ \hline
    \end{tabular}\\
    \textit{Note}: The estimated ESS does not accurately reflect the true uncertainty because the chains have not yet converged. As a result, these ESS estimates can give a misleading impression of having a large number of effective samples when, in reality, the number of effective samples is small.

\end{table}

\section{Limitations and Future Work}

Despite the acceleration achieved through high-rank Hessian approximations, seismic UQ remains computationally costly. More accurate approximations are needed to improve sampling efficiency. For instance, a greater computational budget could be allocated to the Hessian approximation stage, optimizing the balance between approximation and sampling efforts. 
We note that the proposed methods rely on the assumption that the Hessian exhibits locality in both spatial and frequency domains, which may be violated in models with sharp parameter variations, leading to approximation errors. Moreover, the current approximation methods do not address the low-frequency components of the Hessian or transmission effects, and hence require low-rank correction for UQ applications. Future work could approximate the low- and high-frequency parts of the Hessian at the same time, thereby reducing errors caused by the high-pass filter.

\section{Conclusion}

In this work, we implemented high-rank Hessian approximation methods to accelerate both seismic inversion and uncertainty quantification. 
In addition, we proposed a novel Hessian approximation technique that unifies the well-established PSF and PDO methods into a coherent framework, based on their underlying duality. These approximations specifically target the mid-frequency components of the Hessian, which are particularly challenging for Bayesian seismic inversion. While low-frequency components are typically well-informed (from a good starting model) and easily captured, and high-frequency components are poorly informed but can be controlled through priors, the mid-frequency components are only moderately informed and not fully governed by prior knowledge.

Numerical results demonstrate significant improvements in seismic inversion, including simultaneous recovery across the entire model domain and earlier recovery of deep regions, which are often of primary interest in practical applications. More importantly, the Hessian approximation enhances UQ by enabling the generation of more effective samples using the MCMC-gpCN algorithm. Without a suitable Hessian approximation, samples tend to concentrate in limited regions, underestimating variance and leading to unreliable uncertainty estimates, potentially misleading decision-makers by understating the true level of uncertainty. In contrast, the variance of the samples generated using the high-rank Hessian approximations more accurately reflects the variance of the posterior distribution.
The proposed methods are also extensible to 3D problems, where the $\Psi$DO framework remains valid and efficient techniques for computing symbol rows and columns continue to apply.

\begin{acknowledgments}
    We extend our sincere gratitude to TotalEnergies for sponsoring this project. Additional support was provided by DOE ASCR grant DE-SC0023171.
    %We also wish to thank the anonymous reviewers for their interest in our work and their valuable contributions. % Include after we get reviews back
\end{acknowledgments}

\section*{AUTHOR CONTRIBUTIONS}
Ghattas and Nammour identified the overarching research need and provided overall supervision.  
Nammour implemented the wave simulation code and developed the methodology for the pseudo-differential operator probing (PDO) method (Section \ref{sec:pdo}).  
Hu and Alger collaboratively implemented the point spread function (PSF) method.  
Hu adapted the PDO method, developed the methodology for the novel low-rank symbol approach, and carried out its implementation.  
Hu and Alger jointly developed the preconditioned inversion solver.  
Hu conducted the uncertainty quantification analysis and performed the result interpretation.  
Hu drafted the original manuscript; all authors contributed to the review and editing of the final version.

\begin{dataavailability}
The Marmousi model used in our experiments is open source and can be downloaded, for example, from the \href{https://wiki.seg.org/wiki/Open_data}{SEG Open Data repository}. The source code for our implementation of the Hessian approximation methods (PSF, PDO, and PSF+), the preconditioned L-BFGS solver, and the MCMC-gpCN sampler is available on GitHub at \href{https://github.com/mathewgaohu/lrs-psido}{\texttt{github.com/mathewgaohu/lrs-psido}}. Results related to the presented tables and figures are included. However, the full solution history of the L-BFGS iterations and the sample history of the MCMC chains are not provided due to large file sizes, but they can be reproduced by running the source code.

The wave simulation code used in our experiments is proprietary to TotalEnergies and is therefore not publicly available. Nonetheless, it can be substituted with any standard wave simulation code that supports computation of the misfit gradient and Hessian-vector products.
\end{dataavailability}

\bibliographystyle{gji}
\bibliography{references}

\appendix{}

\section{Additional details of the pseudo-differential operator probing method (PDO)}

\subsection{Selection of sample frequencies $\vec \xi_k = (\rho_0, \theta_k)$ on a circle}

Since the seismic source is band-limited, the misfit Hessian contains significant information only within this frequency band. To ensure meaningful symbol approximation, the sampled frequencies should lie within this effective band. Suppose the source function is band-limited to the interval $[f_{\min}, f_{\max}]$. Using the relationship $c = f \lambda$, the corresponding valid range of the sample frequency magnitude $\rho_0$ at a spatial point $\vec x$ is approximately
\begin{equation}
    \frac{f_{\min}}{c(\vec x)} \leq \rho_0 \leq \frac{f_{\max}}{c(\vec x)}.
\end{equation}
To satisfy this condition for all spatial points, $\rho_0$ should lie within the range:
\begin{equation}
    \frac{f_{\min}}{c_{\min}} \leq \rho_0 \leq \frac{f_{\max}}{c_{\max}}.
\end{equation}
Note that this imposes a requirement for a large enough frequency band:  $f_{\max}/f_{\min} > c_{\max}/c_{\min}$.

Accurate angular interpolation requires a sufficient number of sampling angles. However, to prevent overlap in the spectrum of $\mathcal{H}_d v_p$, the sample frequencies must be well-separated to ensure that the spectral components of different symbol columns remain distinct. To maximize angular resolution without causing overlap, we typically choose a large $\rho_0 \approx f_{\max}/c_{\max}$ to increase available spacing. Additionally, to construct real-valued probing vectors, we sample angle pairs $(\theta, \theta + \pi)$ so that their corresponding frequency vectors are negatives of each other, and hence the imaginary parts of the Fourier basis cancel out in the probing vector, i.e., $e^{i \vec x \cdot \vec \xi} + e^{-i \vec x \cdot \vec \xi}$ is real. In practice, using eight equally spaced angles (equivalent to four cosine components) offers a good balance between angular resolution and spectral separability.

\subsection{Implementation overview}

The proposed PDO method is implemented as follows:

\begin{enumerate}
    \item Choose frequency samples $\vec \xi_k = (\rho_0, \theta_k)$ based on the source frequency band. Verify that the frequency band condition $f_{\max}/f_{\min} > c_{\max}/c_{\min}$ is satisfied, and set $\rho_0 \approx f_{\max}/c_{\max}$. Use $r = 8$ equally spaced angles $\theta_k = 2\pi k / r$, for $k=1,2,\ldots,r$.
    \item Construct the probing vector $v_p$ as defined in Equation~(\ref{equ:pdo-vp}), apply $\mathcal{H}_d$ to $v_p$, and extract the symbol columns via spectral separation and inverse FFT.
    \item Construct the interpolation weights $w_k(\vec \xi)$ according to Equation~(\ref{equ:pdo-w}).
    \item Assemble the low-rank symbol using Equation~(\ref{equ:pdo-symbol}) and apply the approximate $\mathcal{H}_d$ using the fast implementation described in Equation~(\ref{equ:lr-psido}), with computational cost $\mathcal{O}(rN\log N)$.
\end{enumerate}

\subsection{Limitations on high-frequency components}
The $\Psi$DO symbol is assumed to grow linearly with frequency. However, due to discretization effects and the band-limited nature of the seismic source, the actual symbol decays in the high-frequency region \citep{demanet2012matrix}. The PDO method tends to overestimate these high-frequency components. While this is generally harmless for gradient preconditioning (since the gradient itself is band-limited), it poses challenges for uncertainty quantification. In particular, computing samples, from the Laplace approximation~(\ref{equ:la-post}), of the form $(\widetilde{H}_d + R)^{-\frac{1}{2}} \vec{z}$ becomes problematic, as the random vector $\vec{z} \sim \mathcal{N}(0,I)$ is not band-limited, leading to inaccuracies in the resulting samples.

\bsp % ``This paper has been produced using the Blackwell
     %   Publishing GJI \LaTeXe\ class file.''

\label{lastpage}

\end{document}